%% file: ms.tex
\newif\ifarxiv
\arxivtrue
\newif\ifprintformalia
\printformaliafalse
\ifarxiv\printformaliafalse\fi

\ifarxiv\else
\fi

\ifarxiv
	\documentclass[parskip=half,10pt,leqno]{amsart}
\else
	\documentclass[parskip=half,12pt,leqno,twoside]{scrreprt}
\fi

\ifarxiv
	\usepackage[hidelinks]{hyperref}
\else
	\usepackage[a-1b]{pdfx}
	\hypersetup{hidelinks}
\fi

\input{other/BA-packages.tex}
\input{other/BA-thesis.tex}
\input{other/BA-math.tex}

\ifarxiv
\fi
\ifarxiv
	\title{A sieve formula for chains of $p$-subgroups}
\else
	\title{\boldmath A sieve formula \\ for chains of $p$-subgroups}
	\subtitle{\boldmath Extending Wielandt's proof of Sylow-Frobenius \\ to a congruence modulo $p^{l+1}$}
\fi

\ifarxiv
	\author{\href{mailto:st173797@stud.uni-stuttgart.de}{Elias Schwesig}}
\else
	\author{Elias Schwesig}
\fi

\ifarxiv
	\date{July 18, 2024}
\else
	\date{July 2024}
\fi

\begin{document}
	\ifarxiv\else
		\pagenumbering{roman}
	\fi
	
	\selectlanguage{english}
	
	\ifarxiv
		\maketitle
		
		\let\thefootnote\relax
		
		\begin{abstract}
			Given a finite group $G$ and a prime $p$, we establish the sieve formula, which is a congruence containing as summands numbers of chains of $p$-subgroups of $G$ of certain orders. This generalises the Theorem of Sylow-Frobenius, using Wielandt's approach. Its name stems from the sieve formula from set theory because of formal similarities.
		\end{abstract}
	\else
		\thispagestyle{empty}
		\vspace*{1cm}
		\begin{center}
			\Huge\bfseries
			\makeatletter
			\@title
			\makeatother
		\end{center}
		\vspace{-2mm}
		\begin{center}
			\Large\bfseries
			\makeatletter
			\@subtitle
			\makeatother
		\end{center}
%		\begin{center}
%			\Huge
%			\color{gray}
%			$$\sum_{\substack{I \subseteq [0, l] \\ I \neq \emptyset}} (-1)^{|I| + 1} \Num(s - I) \equiv_{p^{l+1}} 1$$
%		\end{center}
		\vfill
		\begin{center}
			\Large
			\makeatletter
			Bachelor thesis
			\makeatother
		\end{center}
		\vfill
		\begin{center}
			\Large
			Elias Schwesig
			
			July 2024
			
			\vfill
			
			\normalsize
%			Submitted to the
%			\vspace{15pt}
			
			% Vorliegende Datei muss als Farbraum RGB haben! Das originale Unilogo liegt in CYMK vor.
			% Das kann man mittels Adobe Acrobat Pro konvertieren ("Convert Colors to Output Intent Adobe RGB (1998)")
			\includegraphics[trim=0cm 2cm 0cm 2cm]{../unilogo-englisch-breit.pdf}
		\end{center}
		\normalsize
		\newpage
		\ \newpage
	\fi
	
	\tableofcontents
	
	%%%%%%%%%%%%%%%%%%%%%%%%%%%%%%%%%%%%%%%%%%%%%%%%%%%%%%%%%%%%%%%%%%%%%%%%%%%%%%%%%%%%%%%
	
	\ifarxiv
		\section*{Conventions}
		\input{text/BA-Conventions.tex}
		
%		\section{Introduction}
		\pagenumbering{arabic}
		\input{text/BA-Introduction.tex}
		
		\section{Preliminaries}
		\input{text/BA-Preliminaries.tex}
		
		\section{A decomposition into orbits}	
		\input{text/BA-A-decomposition-into-orbits.tex} % Hier keine Leerzeile lassen!
		\subsection{Subsets, $p$-subgroups and chains of $p$-subgroups}
		\input{text/BA-subsets-p-subgroups.tex}
		
		\subsection{Orbit lengths}
		\input{text/BA-orbit-lengths.tex}
		
		\subsection{Wielandt's lemma}
		\input{text/BA-Wielandt-lemma.tex}
		
		\subsection{Sylow-Frobenius}
		\input{text/BA-sub-formula-a_0.tex}
		
		\section{Counting orbit representatives}
		Let $s \in [0, t]$.
		
		\subsection{Transversals}
		\input{text/BA-Transversals.tex}
		
		\subsection{Constructions}
		\input{text/BA-Omega-phi-construction}
		
		\subsection{Some properties of the assembly map}
		\input{text/BA-Some_basic_properties.tex}
		
		\subsection{A formula for \texorpdfstring{$a_l$}{a\_l}}
		\input{text/BA-a_l.tex}
		
		\section{Counting $p$-subgroups}
		\subsection{The sieve formula}
		\input{text/BA-sub-formula-a_l.tex}
		
		\subsection{A shortcut using Sylow?}
		\input{text/BA-coherence-Sylow.tex}
		
		\section{Examples}
%		\begin{note}
%			In the following, we will discuss \autoref{thm-sub-formel-pl} for certain values of $l$. We will give direct proofs for $l \in \{1, 2\}$, which is possible without induction and which might be useful for one who is interested in these particular cases.
%		\end{note}
		
%		\subsection{Case $l = 0$}
%		\input{text/BA-a_0.tex}
		
		\subsection{Case $l = 1$}
		\input{text/BA-sub-formula-a_1.tex}
		
		\subsection{Case $l = 2$}
		\input{text/BA-sub-formula-a_2.tex}
		
		\section{Further questions}
		\input{text/BA-questions.tex}

\input{text/BA-Bibliography.tex}
		\ifprintformalia
		\input{text/BA-formalia.tex}
		\fi
		
	%%%%%%%%%%%%%%%%%%%%%%%%%%%%%%%%%%%%%%%%%%%%%%%%%%%%%%%%%%%%%%%%%%%%%%%%%%%%%%%%%%%%%%%
	
	\fi
	%%%%%%%%%%%%%%%%%%%%%%%%%%%%%%%%%%%%%%%%%%%%%%%%%%%%%%%%%%%%%%%%%%%%%%%%%%%%%%%%%%%%%%%
\end{document}

%% file: other/BA-packages.tex
\usepackage[utf8]{inputenc}
\usepackage{etoolbox}

\usepackage{minibox} % für Boxen, die auch mehrere Zeilen enthalten dürfen

\makeatletter
\usepackage{enumitem}

\makeatother

\ifarxiv\else
	\usepackage{scrlayer-scrpage}
	\makeatletter
	\addtokomafont{disposition}{\rmfamily} % Standard-Schriftart in chapters, sections etc.
	\makeatother
\fi

\usepackage{scrhack}

\usepackage{xspace}

\usepackage{xparse}

\usepackage{relsize}
\usepackage{exscale}

\usepackage{amsmath}
\usepackage{amssymb}
\usepackage{amsthm}
\usepackage{mathtools}

\ifarxiv
\else
	\usepackage{lmodern}
	\usepackage[T1]{fontenc}

\fi

\makeatletter
\let\oldtheHtheorem\theHtheorem
\renewcommand{\theHtheorem}{\textup{\oldtheHtheorem}}
\makeatother

\usepackage{thmtools}
\usepackage{thm-autoref}

\usepackage{wasysym}
\usepackage{stmaryrd}

\usepackage{tikz}
\usetikzlibrary{arrows, automata, positioning, graphs, decorations.pathreplacing, decorations.pathmorphing, matrix, babel}

\usepackage{tikz-cd} % für dieses Paket wird babel benötigt

\AtBeginEnvironment{align}{\setcounter{equation}{0}} % Zählungen bei align werden immer neu begonnen

\usepackage{eqparbox}

\usepackage[textheight=23cm, % ursprünglich: 22.6cm
			textwidth=16cm,
			hcentering]{geometry}
			
\usepackage[hang,flushmargin]{footmisc}

%% file: other/BA-thesis.tex
\usepackage[english,german]{babel}
\usepackage{hyphenat}
\usepackage[noadjust]{cite}

\makeatletter
\newcommand{\bfcite}[2]{{\bfseries\cite[\mdseries #2]{#1}}}
\newcommand{\bfcit}[1]{{\bfseries\cite{#1}}}
\makeatother

\declaretheoremstyle[bodyfont=\upshape,headpunct={},postheadspace=10pt]{upright}
\declaretheoremstyle[bodyfont=\itshape,headpunct={},postheadspace=10pt]{italic}

\declaretheorem[style=upright,name=Definition\upshape]{definition}
\declaretheorem[style=upright,sibling=definition,name=Question\upshape]{question}
\declaretheorem[style=upright,sibling=definition,name=Reminder\upshape]{reminder}
\declaretheorem[style=upright,sibling=definition,name=Notation\upshape]{notation}
\declaretheorem[style=upright,sibling=definition,name=Remark]{remark without proof}
\declaretheorem[style=upright,sibling=definition,name=Example\upshape]{example}
\declaretheorem[style=italic,sibling=definition,name=Theorem\upshape]{theorem}
\declaretheorem[style=italic,sibling=definition,name=Lemma\upshape]{lemma}
\declaretheorem[style=italic,sibling=definition,name=Remark\upshape]{remark}
\declaretheorem[style=italic,sibling=definition,name=Corollary\upshape]{corollary}
\declaretheorem[style=italic,sibling=definition,name=Proposition\upshape]{proposition}

%% file: other/BA-math.tex
\newcommand{\ie}{\text{i.e.}\xspace}
\newcommand{\cf}{\text{cf.}\xspace}
\newcommand{\Cf}{\text{Cf.}\xspace}

\newcommand{\Ad}{\textit{Ad}\xspace}

\makeatletter
\renewcommand{\leq}{\leqslant}
\renewcommand{\geq}{\geqslant}
\renewcommand{\nleq}{\nleqslant}

\newenvironment{note}
	{\begin{quote} \footnotesize\upshape}
	{

	\end{quote}}
\newcommand{\mystackrel}[3][T]{\stackrel{\eqmakebox[#1]{\scriptsize#2}}{#3}}

\makeatletter
\NewDocumentCommand{\refstackrel}{O{} O{=} m}{%
	\edef\@theoremname{\expandafter\@car\string#3\@nil}%
	\protected@edef\@currentlabelname{\MakeUppercase{\@theoremname}}%
	\stackrel{\text{\@currentlabelname.\,\ref{#3}\xspace#1}}{#2}%
}
\makeatother

\mathcode`l="8000
\begingroup
\makeatletter
\lccode`\~=`\l
\DeclareMathSymbol{\lsb@l}{\mathalpha}{letters}{`l}
\lowercase{\gdef~{\ifnum\the\mathgroup=\m@ne \ell \else \lsb@l \fi}}%
\endgroup

\newcommand{\textover}[3][l]{%
	\makebox[\widthof{$#2$}][#1]{$#3$}%
}

\makeatletter
\def\qed@tag@alignat{%
	\global\tag@true \nonumber
	&\omit\setboxz@h {\strut@ \qedsymbol}%
	\iftagsleft@%
	\global\advance\tagshift@-\displaywidth%
	\fi%
	\tagsleft@false
	\place@tag
	\kern-\tabskip
	\ifst@rred \else \global\@eqnswtrue \fi \global\advance\row@\@ne \cr
}

\def\alignat@qed{%
	\ifmeasuring@ \tag*{\qedsymbol}%
	\else \let\math@cr@@@\qed@tag@alignat
	\fi
}
\@xp\let\csname alignat*@qed\endcsname\alignat@qed
\makeatother

\newcommand{\Z}{\mathbb Z}
\newcommand{\Q}{\mathbb Q}
\newcommand{\F}{\mathbb F}

\newcommand{\val}{\operatorname{v}}

\newcommand{\Stab}{\operatorname{Stab}}

\NewDocumentCommand{\spn}{O{} m}{%
	\langle #2 \IfValueT{#1}{\rangle_{#1}}%
}

\newcommand{\id}{\operatorname{id}}

\newcommand{\Sym}{\operatorname{S}}

\newcommand{\Alt}{\operatorname{A}}

\newcommand{\Cycl}{\mathrm{C}}

\newcommand{\GL}{\operatorname{GL}}

\newcommand{\SL}{\operatorname{SL}}

\newcommand{\C}{\mathrm C}

\newcommand{\Sub}{\operatorname{Sub}}
\newcommand{\Num}{\operatorname{N}}
\newcommand{\Tv}{\operatorname{Tv}}
\newcommand{\tv}{\operatorname{tv}}
\newcommand{\Data}{\Theta}
\newcommand{\TvBin}{\mathrm b}

\usepackage{upgreek}
\newcommand{\DataPhi}{\upvarphi}
\newcommand{\DataRho}{\uprho}

%% file: text/BA-Conventions.tex
Let $X, Y$ be sets. Let $G$ be a group. Let $R$ be a commutative ring.
\begin{itemize}[listparindent=0pt]\ifarxiv\else\small\fi
	\item We write $[a, b] := \{z \in \Z: a \leq z \leq b\}$ for $a, b \in \Z$.
	\item Given $n \in \Z$, we write $\Z_{\geq n} := \{z \in \Z: z \geq n\}$.
	\item Given $x \in \Z$ and $T \subseteq [a, b]$ for some $a, b \in \Z$, we write $x - T := \{x - t : t \in T\}$. For example, we have $3 - \{2, 6\} = \{1, -3\}$ and $5 - \emptyset = \emptyset$.
	\item Given $a, b, x \in R$, we write $a \equiv_x b$ to indicate that $a - b \in xR$.
	\item Given $n \in \Z$, we write $\val_p(n) := \max\{k \in \Z_{\geq 0}: n \equiv_{p^k} 0\}$ if $n \neq 0$, and $\val_p(0) := \infty$.
	
	Let $z + \infty := \infty$ for $z \in \Z$ and $\infty + \infty := \infty$.
	
	Furthermore, we stipulate that $\infty \geq z$ for $z \in \Z$.
	
	\item We write $\mathfrak P(X) := \{Y: Y \subseteq X\}$ for the power set of $X$.
	\item We write $X \subseteq Y$ if $X$ is a subset of $Y$.
	\item We write $X \sqcup Y := X \cup Y$ if the union is disjoint, \ie $X \cap Y = \emptyset$. Furthermore, given a set $U$ and a tuple $(A_i)_{i \in I}$ of subsets $A_i \subseteq U$ for $i \in I$ which are pairwise disjoint, then we write $$\bigsqcup_{i \in I} A_i := \bigcup_{i \in I} A_i.$$
	\item Let $(A_i)_{i \in I}$ be a tuple of sets. Then we write $$\coprod_{i \in I} A_i := \bigsqcup_{i \in I} \{(i, a): a \in A_i\} = \{(i, a): i \in I, a \in A_i\}.$$
	\item We write $H \leq G$ is $H$ is a subgroup of $G$, and $H < G$ if $H \leq G$ and $H \neq G$.
	\item Let $M \subseteq G$. In this case, we write $M \nleq G$ to indicate that $M$ is a subset, but not a subgroup of $G$.
	\item Given $g, x \in G$, we write ${}^xg := xgx^{-1}$ and $g^x := x^{-1}gx$.
	
	Given $g \in G$ and $H \leq G$, we write ${}^gH := \{{}^gh : h \in H\}$ and $H^g := \{h^g : h \in H\}$.
	
	\item Let $X$ be a $G$-set. Let $x \in X$. We write $\Stab(x) := \Stab_G(x) := \{g \in G: g\cdot x = x\}$ for the stabilizer of $x$ in $G$.
%	\item Let $G$ be finite. Let $p$ be a prime. Write $|G| = p^tn$, where $t \in \Z_{\geq 0}$, $n \in \Z_{\geq 1}$ and $n \not \equiv_p 0$. Then we write $\Syl_p(G) = \{H \leq G: |H| = p^t\}$.
	\item Let $n \in \Z_{\geq 1}$. In the symmetric group $\Sym_n$, we write $\sigma \cdot \tau \in \Sym_n$ for the composite of $\sigma \in \Sym_n$ and $\tau \in \Sym_n$, where first $\sigma$ and then $\tau$ is applied, e.g. $(1, 2) \cdot (2, 3) = (1, 3, 2)$ in $\Sym_3$.\begin{note}
		We do this to follow the convention used in Magma \bfcit{Mag}.
	\end{note}
\end{itemize}

%% file: text/BA-Introduction.tex
\ifarxiv
	\section{The main result}
\fi
Let $G$ be a finite group. Let $p$ be a prime.

We write $|G| = p^tn$, where $t = \val_p(|G|) \in \Z_{\geq 0}$, $n \in \Z_{\geq 1}$ and $n \not \equiv_p 0$.

Let $s \in [0, t]$.

\ifarxiv\else
	\section{The main result}\label{sec-intro-main-result}
\fi
Let $k \geq 0$. Let $I = \{c_1, \dots, c_k\} \subseteq [0, t]$, where $c_1 > \ldots > c_k$. We write $$\Num(s-I) := |\{(U_1, \dots, U_k) : U_1 \leq U_2 \leq \ldots \leq U_k \leq G,\ |U_i| = p^{s-c_i} \text{ for } i \in [1,k]\}|.$$ That is, $\Num(s-I)$ denotes the number of chains of $p$-subgroups $$U_1 \leq U_2 \leq \ldots \leq U_k \leq G,$$ where $U_i$ has order $p^{s-c_i}$ for $i \in [1,k]$.

Then we obtain the following sieve formula for $l \in [0, s]$\ifarxiv; \cf \autoref{thm-sub-formel-pl}\fi.
\begin{align*}
	\sum_{\substack{I \subseteq [0, l] \\ I \neq \emptyset}} (-1)^{|I| + 1} \Num(s-I) \equiv_{p^{l+1}} 1.
\end{align*}

\ifarxiv
We have used the approach Wielandt employed to give a short proof of the Theorem of Sylow-Frobenius \bfcit{Wie}. \Cf also \autoref{remark-q-mod-p^k}.
\fi

Its name stems from the sieve formula from set theory, also known as inclusion-exclusion principle, because of formal similarities.

%% file: text/BA-Preliminaries.tex
\ifarxiv\else
	\section{Disjoint subsets}
\fi
\begin{remark}\label{remark-disjoint-union-index-sets}
	Let $X$, $I$ be sets. Let $X_i \subseteq X$ be non-empty subsets for $i \in I$.
	
	Suppose that $X = \bigsqcup_{i \in I} X_i$, \ie $X = \bigcup_{i \in I} X_i$ and $X_i \cap X_j = \emptyset$ for $i,j \in I$ with $i \neq j$.
	
	Let $J, \hat J \subseteq I$ with $\bigsqcup_{j \in J} X_j = \bigsqcup_{j \in \hat J} X_j$. Then we have $J = \hat J$.
	
	\begin{proof}
		We show $\stackrel !\subseteq$. The other inclusion follows vice versa.
		
		Suppose given $j \in J$. Choose $x \in X_j$, which is possible since $X_j \neq \emptyset$.
		
		Since $X_j \subseteq \bigsqcup_{j \in J} X_j = \bigsqcup_{j \in \hat J} X_j$, we have $x \in \bigsqcup_{j \in \hat J} X_j$. So we may choose $k \in \hat J$ such that $x \in X_k$.
		
		So $x \in X_j \cap X_k$. This yields $X_j \cap X_k \neq \emptyset$. Hence $j = k \in \hat J$.
	\end{proof}
\end{remark}

\ifarxiv\else
	\section{Congruences}
\fi
\begin{definition}
	Let $p$ be a prime. Let $$\Z_{(p)} := \{\tfrac ab: a \in \Z,\ b \in \Z \setminus p\Z\} \subseteq \Q.$$ Note that $\Z \subseteq \Z_{(p)}$.
\end{definition}

\begin{remark}\label{rem-z-p-ringiso}
	Suppose given a prime $p$ and $k \in \Z_{\geq 0}$.
	\begin{enumerate}[listparindent=0pt]
		\item Let $a \in \Z$ and $b \in \Z\setminus p\Z$. Then we have $\frac ab \in p^k\Z_{(p)}$ if and only if $a \in p^k\Z$.
		\item The following ring morphism is an isomorphism.
		\begin{align*}
			\Z/p^k\Z &\to \Z_{(p)}/p^k\Z_{(p)}, \\
			z + p^k\Z &\mapsto \tfrac z1 + p^k\Z_{(p)}
		\end{align*}
		In particular, if $a, b \in \Z$ are given such that $a \equiv_{p^k} b$ holds in $\Z_{(p)}$, then $a \equiv_{p^k} b$ holds in $\Z$.
	\end{enumerate}
	\begin{proof}
		\Ad (1). Suppose that $a \in p^k\Z$. We may choose $z \in \Z$ with $a = p^kz$. Then $$\tfrac ab = \tfrac{p^kz}{b} = p^k \tfrac zb \in p^k\Z_{(p)}.$$
		
		Conversely, suppose that $\frac ab \in p^k\Z_{(p)}$. We may choose $c \in \Z$ and $d \in \Z \setminus p\Z$ with $\frac ab = p^k \frac cd$. It follows that $ad = p^kcb$ and so $$\val_p(a) = \val_p(a) + \val_p(d) = \val_p(ad) = \val_p(p^kcb) \geq k.$$ Therefore, $a \in p^k\Z$.
		
		\Ad (2). \emph{Injective}. Suppose given $z \in \Z$ with $z + p^k\Z_{(p)} = 0$. Then $\frac z1 \in p^k\Z_{(p)}$ and, by (1), $z \in p^k\Z$, \ie $z + p^k\Z = 0$.
		
		\emph{Surjective}. Suppose given $\frac ab \in p^k\Z_{(p)}$ with $a \in \Z$ and $b \in \Z \setminus p\Z$. Then $\gcd(b,p^k) = 1$. Bezout's Lemma gives $s, t \in \Z$ such that $sb + tp^k = 1$. Then $sb - 1 = tp^k$. Hence $$a(sb - 1) = atp^k \in p^k\Z.$$
		
		By (1), this is equivalent to $\frac{a(sb - 1)}{b} = as - \frac ab \in p^k\Z_{(p)}$.
		Hence \[as + p^k\Z \ \mapsto\ as + p^k\Z_{(p)} = \tfrac ab + p^k\Z_{(p)}.\qedhere\]
		%			Since $z-w \in \Z$, we have $z-w \in p^k\Z_{(p)} \cap \Z$. We have to show $z-w \stackrel !\in p^k\Z$.
		%			
		%			Write $x := z - w$. Then $x = p^k\frac ab$ for some $a \in \Z$, $b \in \Z \setminus p\Z$, \ie $$xb = p^ka \equiv_{p^k} 0.$$ So $$\val_p(xb) \geq k.$$ Since $\val_p(xb) = \val_p(x) + \val_p(b)$ and $\val_p(b) = 0$ since $b \notin p\Z$, it follows that $\val_p(x) \geq k$, \ie $x \in p^k\Z$.
	\end{proof}
\end{remark}

%\begin{lemma}\label{lemma-p-sum}
%	Suppose given a prime $p$. We obtain the following congruence in $\Z_{(p)}$.
%	$$\sum_{i \in [1,p-1]} \tfrac 1i \equiv_p \sum_{i \in [1,p-1]} i.$$
%	
%	\begin{proof}
%		We write \begin{align*}
%			\psi:\quad \Z/p\Z &\to \Z_{(p)}/p\Z_{(p)} \\
%			a + p\Z &\mapsto \psi(a + p\Z) := \tfrac a1 + p\Z_{(p)}
%		\end{align*} for the ring isomorphism from \autoref{rem-z-p-ringiso}.(2).
%		
%		We obtain the following.
%		\begin{align*}\textstyle
%			(\sum_{i \in [1,p-1]} \frac 1i) + p\Z_{(p)}
%			&= \textstyle\sum_{i \in [1,p-1]} (\frac 1i + p\Z_{(p)}) \\
%			&= \textstyle\sum_{i \in [1,p-1]} (i + p\Z_{(p)})^{-1} \\
%			&= \textstyle\sum_{i \in [1,p-1]} (\psi(i + p\Z))^{-1} \\
%			&= \textstyle\psi(\sum_{i \in [1,p-1]} (i + p\Z)^{-1}) \\
%			&= \textstyle\psi(\sum_{x \in (\Z/p\Z)^\times} x^{-1}) \\
%			&= \textstyle\psi(\sum_{x \in (\Z/p\Z)^\times} x) \\
%			&= \textstyle\psi(\sum_{i \in [1,p-1]} (i + p\Z)) \\
%			&= \textstyle\sum_{i \in [1,p-1]} \psi(i + p\Z) \\
%			&= \textstyle\sum_{i \in [1,p-1]} (i + p\Z_{(p)}) \\
%			&= \textstyle(\sum_{i \in [1,p-1]} i) + p\Z_{(p)}. \qedhere
%		\end{align*}
%	\end{proof}
%\end{lemma}

\ifarxiv\else
	\section{Noetherian induction}
\fi
Let $(X, \leq)$ be a partially ordered set.

\begin{definition}\leavevmode
	\begin{enumerate}[listparindent=0pt]
		\item Let $T \subseteq X$. We say that $t \in T$ is \emph{minimal} in $T$ if $\{s \in T: s < t\} = \emptyset$.
		\item We say that $(X, \leq)$ is \emph{noetherian} if every non-empty subset $T \subseteq X$ has a minimal element in $T$.
	\end{enumerate}
\end{definition}

\begin{lemma}\label{lemma-noetherian-chains} The following assertions are equivalent.
	\begin{enumerate}[listparindent=0pt]
		\item The partially ordered set $(X, \leq)$ is \textsf{not} noetherian.
		\item There exist $x_n \in X$ for $n \in \Z_{\geq 1}$ such that $x_i > x_{i+1}$ for $i \in \Z_{\geq 1}$.
	\end{enumerate}
	
	\begin{proof}
		\Ad (1) $\Rightarrow$ (2). Since $(X, \leq)$ is not noetherian, we may choose a non-empty subset $T \subseteq X$ with no minimal element.
		
		Choose $x_1 \in T$. Then $x_1$ is not minimal. Hence, we may choose $x_2 \in T$ such that $x_1 > x_2$.
		
		Then $x_2$ is not minimal. Hence, we may choose $x_3 \in T$ such that $x_2 > x_3$.
		
		Repeating this, we get a infinite decreasing sequence $$x_1 > x_2 > x_3 > \dots$$ as desired.
		
		\Ad (2) $\Rightarrow$ (1). Choose $x_n \in X$ for $n \in \Z_{\geq 1}$ such that $x_i > x_{i+1}$ for $i \in \Z_{\geq 1}$. Then the non-empty subset $T := \{x_n : n \in \Z_{\geq 1}\} \subseteq X$ has no minimal element in $T$, hence $(X, \leq)$ is not noetherian.
	\end{proof}
\end{lemma}

\begin{corollary}\label{cor-finite-posets-noetherian}
	If $X$ is a finite set, then $(X, \leq)$ is noetherian.
	
	\begin{proof}
		Suppose that $(X, \leq)$ is not noetherian. Then we may choose $x_n \in X$ for $n \in \Z_{\geq 1}$ such that $$x_1 > x_2 > x_3 > \dots,$$ \cf \autoref{lemma-noetherian-chains}.(1 $\Rightarrow$ 2). This implies that $|\{x_n : n \in \Z_{\geq 1}\}| = \infty$ and therefore $|X| = \infty$. 
	\end{proof}
\end{corollary}

\begin{lemma}[Noetherian induction]\label{lemma-noetherian-induction}
	Suppose that $(X, \leq)$ is noetherian.
	
	Suppose given a statement $P(x)$ for $x \in X$. 
	
	Suppose that the following property holds.
	
	\begin{enumerate}[listparindent=0pt]
		\item[\upshape{(N)}] If $x \in X$ is given such that $P(x')$ holds for $x' \in X$ with $x' < x$, then $P(x)$ holds. %If $x \in X$ is given such that $\{x' \in X: x' < x\} \subseteq Y$, then $x \in Y$.
	\end{enumerate}
	
	Then $P(x)$ holds for $x \in X$.

	\begin{proof}
	%	Let $Y := \{x \in X : P(x)\}$. \emph{Assume} that $Y \neq X$. Then $X \setminus Y \neq \emptyset$. Hence $Z := X \setminus Y \subseteq X$ has a minimal element $z \in Z$. Then $z \in X$, and $$\{x' \in X : x' < z\}\subseteq X \setminus Z = Y,$$ so $z \in Y$ by (N), which is a \emph{contradiction}.
		\emph{Assume} that $Z := \{x \in X : \lnot P(x)\} \neq \emptyset$. Hence, $Z$ has a minimal element $z \in Z$. Then $z \in X$, and $P(x')$ holds for $x' \in X$ with $x' < z$ since $z$ is minimal in $Z$. By (N), $P(z)$ holds, which is a \emph{contradiction} to $z \in Z$.
	\end{proof}
\end{lemma}

%% file: text/BA-A-decomposition-into-orbits.tex
\ifarxiv
	Throughout the text, we keep the following data:
	
	\minibox[frame]{
		Let $G$ be a finite group. Let $p$ be a prime.\\
		We write $|G| = p^tn$, where $t = \val_p(|G|) \in \Z_{\geq 0}$, $n \in \Z_{\geq 1}$ and $n \not \equiv_p 0$.
	}
\else
	\begin{picture}(0,0)
		\put(-2,150){
			Throughout the text, we keep the following data:
		}
		\put(-5, 120){
			\minibox[frame]{
				Let $G$ be a finite group. Let $p$ be a prime.\\
				We write $|G| = p^tn$, where $t = \val_p(|G|) \in \Z_{\geq 0}$, $n \in \Z_{\geq 1}$ and $n \not \equiv_p 0$.
			}
		}
	\end{picture}
	
	\vspace*{-20pt}
	Suppose given $s \in [0, t]$.
	\vspace*{-5pt}
\fi

%% file: text/BA-subsets-p-subgroups.tex
\begin{definition}\label{definition-sub-num}\leavevmode
	\begin{enumerate}[listparindent=0pt]
		\item We write $$\Omega := \Omega^{[s]} := \{M \subseteq G: |M| = p^s\} \subseteq \mathfrak P(G).$$ We remark that $|\Omega^{[s]}| = \binom{p^tn}{p^s}$.
		
		\item Let $\ell \in [0,t]$. Let $$\Sub(\ell) := \{U \leq G: |U| = p^\ell\}.$$
		Let $I = \{c_1, \dots, c_k\} \subseteq [0, t]$ where $c_1 < \ldots < c_k$.
		Let
		$$\Sub(c_1, \dots, c_k) := \{(U_1,\dots, U_k) \in \Sub(c_1) \times \ldots \times \Sub(c_k): U_1 \leq \ldots \leq U_k\}.$$
		We also write $$\Sub(I) := \Sub(c_1, \dots, c_k).$$
		Let
		$$\Num(c_1, \dots, c_k) := \Num_{G,p}(c_1, \dots, c_k) := |\Sub(c_1, \dots, c_k)|.$$
		We also write $$\Num(I) := \Num_{G,p}(I) := \Num(c_1, \dots, c_k).$$
		Let $s - I := \{s-i : i \in I\}$. Then $$\Sub(s-I) = \Sub(s-c_k, \dots, s-c_1),$$ and $$\Num(s-I) = \Num(s-c_k, \dots, s-c_1).$$
		Note that $\Sub(\emptyset) = \{()\}$ and so $\Num(\emptyset) = 1$.
		
		\item Let $I = \{c_1, \dots, c_k\} \subseteq [0, t]$ where $c_1 < \ldots < c_k$. Let $V \leq G$. Let
		$$\Sub(I, V) := \{(U_1, \dots, U_k, V) : (U_1, \dots, U_k) \in \Sub(c_1, \dots, c_k),\ U_k \leq V\}.$$
	\end{enumerate}
\end{definition}

\begin{example}\label{ex-sym-4-sub-num}
	Let $G = \Sym_4$ and $p = 2$.
	
	We have the following partially ordered set of $2$-subgroups in $\Sym_4$.
	
	\tiny
	\newcommand{\cm}{{,}}
	\begin{tikzcd}[column sep=-1ex,row sep=15ex]
		&&&&&\Sym_4\arrow[no head, dl]\arrow[no head, d]\arrow[no head, dr] \\
		&&&&\spn{(1\cm2),(1\cm3)(2\cm4)}
		\arrow[no head, dlll]\arrow[no head, d]\arrow[no head, drrr]
		&\spn{(2\cm4),(1\cm2)(3\cm4)}
		\arrow[no head, dlll]\arrow[no head, d]\arrow[no head, drr]
		&\spn{(2\cm3),(1\cm3)(2\cm4)}
		\arrow[no head, dlll]\arrow[no head, d]\arrow[no head, dr] \\
		&\spn{(1\cm2),(3\cm4)}
		\arrow[no head, dl]\arrow[no head, drrrr]\arrow[no head, drrrrr]
		&\spn{(1\cm3),(2\cm4)}
		\arrow[no head, dl]\arrow[no head, drr]\arrow[no head, drrrrr]
		&\spn{(1\cm4),(2\cm3)}
		\arrow[no head, dl]\arrow[no head, d]\arrow[no head, drrrrr]
		&\spn{(1\cm3\cm2\cm4)}
		\arrow[no head, drr]
		&\spn{(1\cm2\cm3\cm4)}
		\arrow[no head, drr]
		&\spn{(1\cm2\cm4\cm3)}
		\arrow[no head, drr]
		&\spn{(1\cm2)(3\cm4), (1\cm4)(2\cm3)}
		\arrow[no head, dl]\arrow[no head, d]\arrow[no head, dr] \\
		\spn{(1\cm2)}
		\arrow[no head, drrrrr]
		&\spn{(1\cm3)}
		\arrow[no head, drrrr]
		&\spn{(1\cm4)}
		\arrow[no head, drrr]
		&\spn{(2\cm3)}
		\arrow[no head, drr]
		&\spn{(2\cm4)}
		\arrow[no head, dr]
		&\spn{(3\cm4)}
		\arrow[no head, d]
		&\spn{(1\cm2)(3\cm4)}
		\arrow[no head, dl]
		&\spn{(1\cm3)(2\cm4)}
		\arrow[no head, dll]
		&\spn{(1\cm4)(2\cm3)}\arrow[no head, dlll] \\
		&&&&&1
	\end{tikzcd}
	
	\normalsize
	\begin{enumerate}[listparindent=0pt]
		\item We have
		\begin{align*}
			\Sub(3) &= \{U \leq \Sym_4 : |U| = 2^3\} \\
			&= \{\spn{(1,2),(1,3)(2,4)},
			\spn{(2,4),(1,2)(3,4)},
			\spn{(2,3),(1,3)(2,4)}\}.
		\end{align*}
		\item We have $\Num(3) = \Num(\{3\}) = |\Sub(3)| = 3$.
		\item We see that $$1 < \spn{(1,2)} < \spn{(1,2),(3,4)}$$ is a chain of $2$-subgroups of order $2^0$, $2^1$, $2^2$ in $\Sym_4$. Hence $$(1, \spn{(1,2)}, \spn{(1,2),(3,4)}) \in \Sub(0, 1, 2) = \Sub([0, 2]).$$
		\item We have $\Num(0,1,2) = |\Sub(0,1,2)| = 15$.
	\end{enumerate}
	
\end{example}

\begin{remark}\label{remark-trivial-subgroup-in-chain}
	Let $I \subseteq [1, t]$. We have $$\Num(I) = \Num(\{0\} \sqcup I).$$
	
	\begin{proof}
		Write $I = \{c_1, \dots, c_k\}$ with $c_1 < \ldots < c_k$.
		
		Note that $\Sub(0) = \{U \leq G : |U| = 1\} = \{1\}$. Therefore
		\begin{align*}
			&\Num(0, c_1, \dots, c_k) \\
			=\ &|\Sub(0, c_1, \dots, c_k)| \\
			=\ &|\{(U_0, U_1,\dots, U_k) \in \Sub(0) \times \Sub(c_1) \times \ldots \times \Sub(c_k): U_0 \leq U_1 \leq \ldots \leq U_k\}| \\
			=\ &|\{(1, U_1,\dots, U_k) \in \{1\} \times \Sub(c_1) \times \ldots \times \Sub(c_k): 1 \leq U_1 \leq \ldots \leq U_k\}| \\
			=\ &|\{(U_1,\dots, U_k) \in \Sub(c_1) \times \ldots \times \Sub(c_k): U_1 \leq \ldots \leq U_k\}| \\
			=\ &\Num(c_1, \dots, c_k).\qedhere
		\end{align*}
	\end{proof}
\end{remark}

\begin{example}
	We have $\Num_{\Sym_4, 2}(1) = \Num_{\Sym_4, 2}(0, 1) = 9$ and $\Num_{\Sym_4, 2}(1, 2) = \Num_{\Sym_4, 2}(0,1,2) = 15$; \cf \autoref{ex-sym-4-sub-num} and \autoref{remark-trivial-subgroup-in-chain}.
\end{example}

%% file: text/BA-orbit-lengths.tex
\begin{definition}
	Let $G$ act on $\Omega$ via multiplication from the left, that is $$g \cdot M := gM = \{gm: m \in M\}$$ for $g \in G$ and $M \in \Omega$.
	
	\begin{enumerate}[listparindent=0pt]
		\item Suppose given $M \in \Omega$. We write $$[M] := \{gM: g \in G\}$$ for the orbit of $M$ under $G$. Furthermore, we write $$\overline \Omega := \overline \Omega^{[s]} := \{[M]: M \in \Omega\}$$ for the set of $G$-orbits.
		\item Let \begin{align*}
			\DataRho: \quad \Omega &\to \overline \Omega \\
			M &\mapsto [M].
		\end{align*}
	\end{enumerate}
\end{definition}

\begin{lemma}\label{lemma-order-stab-and-orbit}
	Suppose given $M \in \Omega$. Let $U := \Stab(M) = \{g \in G: gM = M\} \leq G$.
	\begin{enumerate}[listparindent=0pt]
		\item We have $|[M]| \cdot |U| = p^tn$.
		\item There exists a unique $l \in [0, s]$ such that $|U| = p^{s-l}$.
		
		There exist $m_1, \dots, m_{p^l} \in M$ such that $$M = \bigsqcup_{i \in [1, p^l]} Um_i.$$
		
		We have $|[M]| = p^{t-s+l}n$.
	\end{enumerate}
	
	\begin{proof}
		\Ad (1). This ensues from the orbit-stabilizer theorem since $|G| = p^tn$.
		
		\Ad (2). Let $U$ act on $M$ via multiplication from the left. Given $m \in M$, the orbit $Um$ is a right coset of $U$, which implies $|Um| = |U|$. 
		
		Hence, we get $k \in \Z_{\geq 1}$ and $m_1, \dots, m_k \in M$ such that $M = \bigsqcup_{i \in [1, k]} Um_i$. Thus $$p^s = |M| = \sum_{i \in [1, k]} |Um_i| = \sum_{i \in [1, k]} |U| = k \cdot |U|.$$ So $|U|$ divides $p^s$. Therefore, there is a unique $l \in [0, s]$ with $|U| = p^{s-l}$. This forces $k = p^l$.
		
		By (1), we have $|[M]| = \frac{|G|}{|U|} = \frac{p^tn}{p^{s-l}} = p^{t-s+l}n.$
	\end{proof}
\end{lemma}

\begin{remark}\label{remark-stab-conjugation}
	Let $g \in G$. Let $M \in \Omega$. Then we have $\Stab(gM) = {}^g\Stab(M)$.
	
	In particular, we have $|\Stab(\hat M)| = |\Stab(M)|$ for $\hat M \in [M]$.
	
	\begin{proof}
		We have \begin{align*}
			\Stab(gM)
			&= \{h \in G: hgM = gM\} \\
			&= \{h \in G: g^{-1}hgM = M\} \\
			&= \{h \in G: h^gM = M\} \\
			&= \{h \in G: h^g \in \Stab(M)\} \\
			&= \{h \in G: h \in {}^g\Stab(M)\} \\
			&= {}^g\Stab(M).\qedhere
		\end{align*}
	\end{proof}
\end{remark}

\begin{definition}\label{definition-omega}
	Let $l \in [0, s]$.
	
	\begin{enumerate}[listparindent=0pt]
		\item Let $$\Omega^\ell := \Omega^{[s],l} := \{M \in \Omega^{[s]}: |[M]| = p^{t-s+\ell}n\} \subseteq \Omega.$$
		So $$\Omega = \bigsqcup_{k \in [0, s]} \Omega^k;$$ \cf \autoref{lemma-order-stab-and-orbit}.(2).
		
		Note that $$\Omega^l = \{M \in \Omega : |\Stab(M)| = p^{s-l}\};$$ \cf \autoref{lemma-order-stab-and-orbit}.(1).
		
		\item Let
		\begin{align*}
			\overline \Omega^l &:= \overline \Omega^{[s],l} \\
			&:= \{[M]: M \in \Omega^l\} \\
			&\phantom{:}= \{[M]: M \in \Omega,\ |[M]| = p^{t-s+\ell}n\} \\
			&\phantom{:}= \{[M]: M \in \Omega,\ |\Stab(M)| = p^{s-l}\} \subseteq \overline \Omega.
		\end{align*}
		So $$\overline \Omega = \bigsqcup_{k \in [0, s]} \overline \Omega^k;$$ \cf \autoref{lemma-order-stab-and-orbit}.(2).
		
		Note that for $M \in \Omega^l$ and $\hat M \in [M]$, we have $\hat M \in \Omega^l$.
		
		Furthermore, we write $$\overline \Omega^{[0,l]} := \overline \Omega^{[s],[0,l]} := \bigsqcup_{k \in [0,l]} \overline \Omega^k \subseteq \overline \Omega.$$
		
		\item Let $$a_\ell := |\overline \Omega^\ell| \in \Z_{\geq 0}.$$
	\end{enumerate}
\end{definition}

\begin{remark}\label{remark-omega-sum}
	We have $|\Omega| = p^{t-s}n\sum_{k \in [0,s]} a_kp^k$.
	
	\begin{proof}
		Given $k \in [0, s]$, we may choose representatives $M_{k,1}, \dots, M_{k, a_k} \in \Omega^k$ such that $$\overline \Omega^k = \{[M_{k,j}]: j \in [1, a_k]\}, \qquad \ie \qquad \Omega^k = \bigcup_{M \in \Omega^k} [M] = \bigsqcup_{j \in [1, a_k]} [M_{k, j}].$$
		
		Hence
		$$\Omega = \bigsqcup_{k \in [0,s]} \Omega^k
		= \bigsqcup_{k \in [0,s]} \bigsqcup_{j \in [1, a_k]}[M_{k,j}]$$
		and therefore
		\[|\Omega| = \sum_{k \in [0,s]} \sum_{j \in [1, a_k]} |[M_{k,j}]| = \sum_{k \in [0,s]} \sum_{j \in [1, a_k]} p^{t-s+k}n = \sum_{k \in [0,s]} a_kp^{t-s+k}n = p^{t-s}n \sum_{k \in [0,s]} a_kp^k.\qedhere\]
	\end{proof}
\end{remark}

\begin{definition}\label{definition-q}
	Let $$q := \frac{1}{p^{t-s}n} \binom{p^tn}{p^s} \in \Z_{\geq 1},$$ \cf \autoref{remark-omega-sum}.
\end{definition}

\begin{remark without proof}\label{remark-q-mod-p^k}
	\autoref{remark-omega-sum} yields $$q = \sum_{k \in [0, s]} a_kp^k.$$ So for $l \in [0,s]$, we get $$q \equiv_{p^{l+1}} \sum_{k \in [0,l]} a_kp^k.$$ In particular, we get $$q \equiv_p a_0, \qquad q \equiv_{p^2} a_0 + a_1p, \qquad q \equiv_{p^3} a_0 + a_1p + a_2p^2$$ etc.
	\ifarxiv Wielandt used the congruence $q \equiv_p a_0$ in his proof of the Theorem of Sylow-Frobenius; \cf \autoref{thm-sylow-wielandt}.\fi
\end{remark without proof}

%% file: text/BA-Wielandt-lemma.tex
\begin{note}
	We shall give an interpretation of the right hand side in terms of members of certain chains of $p$-subgroups; \cf \autoref{proposition-al} below.
\end{note}
\begin{remark}\label{remark-stab-of-subgroup}
	Suppose given $H \in \Sub(s) \subseteq \Omega^{[s]}$. Then $\Stab(H) = H$ and $H \in \Omega^0$.
	
	\begin{proof}
		Suppose given $g \in G$. Then $g \in \Stab(H)$ if and only if $gH = H$, \ie $g \in H$. Since $|\Stab(H)| = |H| = p^s$, we get $H \in \Omega^0$, \cf \autoref{definition-omega}.(1).
	\end{proof}
\end{remark}

\begin{lemma}[Wielandt, \cf\ \bfcit{Wie}]\label{lemma-wielandt}
	Let $M \in \Omega$.
	\begin{enumerate}[listparindent=0pt]
		\item Given $M \in \Omega^0$, there is a unique $H \in [M]$ with $H \leq G$.
		\item Given $M \in \Omega^\ell$ for some $l \in [1, s]$, we have $M \nleq G$.
	\end{enumerate}
	
	\begin{proof}
		\Ad (1). We have $|\Stab(M)| = p^{s-0} = p^s$; \cf \autoref{definition-omega}.(1). \autoref{lemma-order-stab-and-orbit}.(2) yields $M = \Stab(M)m$ for some $m \in M$. Hence $m^{-1}M = m^{-1}\Stab(M)m \in [M]$ and $m^{-1}\Stab(M)m = \Stab(M)^m \leq G$.
		
		Suppose given $H, \hat H \in [M]$ such that $H, \hat H \leq G$. We may choose $g \in G$ such that $gH = \hat H$. Since $H = g^{-1}\hat H$, we have $g^{-1} \in H$, thus $g = (g^{-1})^{-1} \in H$. Thus $H = gH = \hat H$.
		
		\Ad (2). \emph{Assume} that $M \leq G$. Then $M \in \Omega^0$, \cf \autoref{remark-stab-of-subgroup}, a \emph{contradiction}.
	\end{proof}
\end{lemma}

\begin{example}
	% Magma-Befehl: G := SymmetricGroup(3); M := {G!1, G!(1, 2)}; Stab(G, M);
	Let $G = \Sym_4$ and $p = 2$.
	
	Let $M = \{\id, (1, 2)\} = \spn{(1, 2)} \in \Omega^{[1]} = \{M \subseteq \Sym_4 : |M| = 2^1\}$.
	
	We have $$\Stab_{\Sym_4}(M) = \{\sigma \in \Sym_4 : \sigma M = M\} = \{\id, (1, 2)\} = M.$$ So $|\Stab_{\Sym_4}(M)| = 2^1 = 2^{1-0}$, hence $M \in \Omega^{[1],0}$.
	
	This also follows by \autoref{remark-stab-of-subgroup}. Note that $$|[M]| = \tfrac{|\Sym_4|}{|\Stab_{\Sym_4}(M)|} = \tfrac{24}{2} = 12.$$ Furthermore, \autoref{lemma-wielandt}.(1) gives $$\{H \leq G: H \in [M]\} = \{M\}.$$
\end{example}

%% file: text/BA-sub-formula-a_0.tex
\begin{note}
	With these tools, Wielandt could already derive \autoref{thm-sylow-wielandt}, which will also result as the particular case $l = 0$ from \autoref{thm-sub-formel-pl}.
\end{note}

\begin{theorem}[\cf\ \bfcite{Syl}{Th.\,II}, \bfcite{Fro}{§4, I.}, \bfcit{Wie}]\label{thm-sylow-wielandt}
	Let $s \in [0, t]$. We have $$\Num(s) \equiv_p 1.$$
	
	\begin{proof}
		Consider $|\Omega| = \binom{p^tn}{p^s} = p^{t-s} \sum_{k \in [0, s]} a_kp^k$, \cf \autoref{remark-omega-sum}. Hence $$q = \tfrac{1}{p^{t-s}} \binom{p^tn}{p^s} = \sum_{k \in [0, s]} a_kp^k \equiv_p a_0,$$ and \autoref{lemma-wielandt}.(1, 2) provides $a_0 = \Num(s)$. That is $q \equiv_p \Num(s)$.
		
		Consider the cyclic group $C := \C_{p^tn}$ of order $p^tn$. It is well-known that $\C_{p^tn}$ has exactly one subgroup of order $p^s$. Therefore \[\Num_{G,p}(s) \equiv_p q \equiv_p \Num_{C,p}(s) = 1.\qedhere\]
	\end{proof}
\end{theorem}

%% file: text/BA-Transversals.tex
Let $l \in [0, s]$.

\begin{remark}\label{remark-omega-1-element-decomposition}
	Let $\hat M \in \Omega^l$. Then there exist $M \in [\hat M]$ and $m_2, \dots, m_{p^l} \in M$ such that $$M = U \sqcup \bigsqcup_{i \in [2,p^l]} Um_i,$$ where $U := \Stab(M) \leq G$.
	
	\begin{proof}
		Thanks to \autoref{lemma-order-stab-and-orbit}.(2), we get $\hat m_1, \dots, \hat m_{p^l} \in \hat M$ such that $\hat M = \bigsqcup_{i \in [1, p^l]} \hat U\hat m_i$ with $\hat U := \Stab(\hat M)$. Let $m_i := \hat m_1^{-1}\hat m_i$ for $i \in [1, p^l]$ and $U := \hat U^{\hat m_1}$. 
		
		Now $M := \hat m_1^{-1} \hat M \in [\hat M]$, and
		\everymath{\textstyle}
		\begin{align*}				
			M &= \hat m_1^{-1}\left(\textstyle\bigsqcup_{i \in [1,p^l]} \hat U\hat m_i\right) \\
			&= \textstyle\bigsqcup_{i \in [1, p^l]} \hat m_1^{-1}\hat U\hat m_i \\
			&= \textstyle\bigsqcup_{i \in [1, p^l]} \hat m_1^{-1}\hat U\hat m_1\hat m_1^{-1}\hat m_i \\
			&= \textstyle\bigsqcup_{i \in [1, p^l]} \hat U^{\hat m_1}\hat m_1^{-1}\hat m_i \\
			&= \hat U^{\hat m_1} \sqcup \textstyle\bigsqcup_{i \in [2, p^l]} \hat U^{m_1}\hat m_1^{-1}\hat m_i \\
			&= U \sqcup \textstyle\bigsqcup_{i \in [2, p^l]} U m_i.
		\end{align*}
		Using \autoref{remark-stab-conjugation}, we get \[U = \hat U^{\hat m_1} = {}^{\hat m_1^{-1}}\hat U = {}^{\hat m_1^{-1}} \Stab(\hat M) = \Stab(\hat m_1^{-1} \hat M) = \Stab(M).\qedhere\]
	\end{proof}
\end{remark}

\begin{reminder}
	Given a subgroup $U \leq G$, we call a subset $T \subseteq G$ a \emph{right transversal} for $U$ in $G$ if $$G = \bigsqcup_{t \in T} Ut.$$
\end{reminder}

\begin{notation}\label{not-tv}
	Given a subgroup $U \leq G$, we choose a right transversal $T \subseteq G$ such that $1 \in T$, and we write $$\tv(U) := \tv_G(U) := T.$$
	
	\begin{note}
		The choice of this transversal, that is the choice of representatives, will affect the calculations, but not the outcome of the considerations.
	\end{note}
	
	We write $$\Tv_{p^l}(U) := \Tv_{G,p^l}(U) := \{T \subseteq (\tv(U) \setminus \{1\}): |T| = p^l - 1\} \subseteq \mathfrak P(\tv(U)).$$ 
\end{notation}

\begin{definition}\label{definition-binomial}
	We write $$\TvBin(l) := \TvBin_s(l) := \binom{p^{t-s+l}n - 1}{p^l - 1} \in \Z_{\geq 0}.$$
\end{definition}

\begin{remark without proof}\label{remark-binomial}
	Let $U \in \Sub(s-l)$ and let $G = \bigsqcup_{t \in \tv(U)} Ut$. Then $$|\tv(U)| = |G/U| = |G|/|U| = p^tn/p^{s-l} = p^{t-s+l}n.$$
	
	We have $|\tv(U)| - 1$ elements in $\tv(U) \setminus \{1\}$.
	
	Any $T \in \Tv_{p^l}(U)$ contains $p^l-1$ elements.
	
	This leads to
	$$\textstyle|\Tv_{p^l}(U)| = \binom{|\tv(U)| - 1}{p^l-1} = \binom{p^{t-s+l}n - 1}{p^l-1} = \TvBin(l).$$
	
	Recall that we have $t = \val_p(|G|)$ and $n = |G|/p^t$.
\end{remark without proof}

\begin{lemma}\label{lemma-b-k-mod-p-k}
	Let $l \in [1, s]$. We have $$\TvBin(l) \equiv_{p^l} \TvBin(l-1).$$
	
	\begin{proof}
		Write $d := t - s \in \Z_{\geq 0}$. We have \allowdisplaybreaks
		\begin{align*}
			\TvBin(l)
			&\mystackrel{D. \ref{definition-binomial}}{=} \binom{p^{t-s+l}n - 1}{p^l - 1} \\
			&\mystackrel{}{=} \frac{\prod_{i \in [1, p^l - 1]} (p^{d+l}n - i)}{\prod_{i \in [1, p^l - 1]} (p^l - i)} \\
			&\mystackrel{}{=} \prod_{i \in [1, p^l - 1]}\frac{p^{d+l}n - i}{p^l - i} \\
			&\mystackrel{}{=} \Biggl(\prod_{\substack{i \in [1, p^l - 1] \\ i \equiv_p 0}} \frac{p^{d+l}n - i}{p^l - i}\Biggr)
			\cdot
			\Biggl(\prod_{\substack{i \in [1, p^l - 1] \\ i \not\equiv_p 0}} \frac{p^{d+l}n - i}{p^l - i}\Biggr) \\
			&\mystackrel{$i=pj$}{=} \Biggl(\prod_{j \in [1, p^{l-1} - 1]}\frac{p^{d+l}n - pj}{p^l - pj}\Biggr)
			\cdot
			\Biggl(\prod_{\substack{i \in [1, p^l - 1] \\ i \not\equiv_p 0}} \frac{p^{d+l}n - i}{p^l - i}\Biggr) \\
			&\mystackrel{}{=} \Biggl(\prod_{j \in [1, p^{l-1} - 1]}\frac{p^{d+(l-1)}n - j}{p^{l-1} - j}\Biggr)
			\cdot
			\Biggl(\prod_{\substack{i \in [1, p^l - 1] \\ i \not\equiv_p 0}} \frac{p^{d+l}n - i}{p^l - i}\Biggr) \\
			&\mystackrel{}{=} \TvBin(l-1)
			\cdot
			\Biggl(\prod_{\substack{i \in [1, p^l - 1] \\ i \not\equiv_p 0}} \frac{p^{d+l}n - i}{p^l - i}\Biggr).
		\end{align*}
		
		Furthermore, we have $$\prod_{\substack{i \in [1, p^l - 1] \\ i \not\equiv_p 0}} \frac{p^{d+l}n - i}{p^l - i}\ \in\ \Z_{(p)}$$ because $p^l - i \not \equiv_p 0$ for $i \in [1, p^l - 1] \setminus p\Z$.
		
		Hence $$\TvBin(l) = \TvBin(l-1) \cdot \Biggl(\prod_{\substack{i \in [1, p^l - 1] \\ i \not\equiv_p 0}} \frac{p^{d+l}n - i}{p^l - i}\Biggr)$$ in $\Z_{(p)}$.
		And $\frac{p^{d+l}n - i}{p^l - i} - 1 = \frac{p^{d+l}n - p^l}{p^l - i} = p^l \frac{p^dn - 1}{p^l - i} \in p^l\Z_{(p)}$ for $i \in [1, p^l - 1] \setminus p\Z$, so $$\frac{p^{d+l}n-i}{p^l-i} \equiv_{p^l} 1$$ for $i \in [1, p^l - 1] \setminus p\Z$, so $$\prod_{\substack{i \in [1, p^l - 1] \\ i \not\equiv_p 0}} \frac{p^{d+l}n-i}{p^l-i} \equiv_{p^l} 1.$$ Hence $$\TvBin(l) = \TvBin(l-1) \cdot \Biggl(\prod_{\substack{i \in [1, p^l - 1] \\ i \not\equiv_p 0}} \frac{p^{d+l}n - i}{p^l - i}\Biggr) \equiv_{p^l} \TvBin(l-1)$$ in $\Z_{(p)}$.
		Thus we obtain the congruence $$\TvBin(l) \equiv_{p^l} \TvBin(l-1)$$ in $\Z$ because $\TvBin(l), \TvBin(l-1) \in \Z$, \cf \autoref{rem-z-p-ringiso}.(2).
	\end{proof}
\end{lemma}

%% file: text/BA-Omega-phi-construction.tex
Let $l \in [0, s]$.

\begin{definition}\label{definition-omega-1-l}\leavevmode
	\begin{enumerate}[listparindent=0pt]
		\item Let
		\begin{align*}
			\Omega_1^l := \Omega_1^{[s],l} &:= \{M \in \Omega^{[s],l}: 1 \in M\} \\
			&\phantom := \{M \in \Omega^{[s]} : 1 \in M,\ |\Stab(M)| = p^{s-l}\};
		\end{align*}
		\cf \autoref{definition-omega}.(1).
		
		\item Let \begin{align*}
			\Omega_1^{[0,l]} &:= \Omega_1^{[s], [0,l]} \\
			&:= \textstyle\{U \sqcup \bigsqcup_{i \in [2,p^l]} Ug_i : U \in \Sub(s-l), \{g_2, \dots, g_{p^l}\} \in \Tv_{p^l}(U)\} \subseteq \Omega^{[s]}.
		\end{align*}
	\end{enumerate}
\end{definition}

\begin{example}\label{ex-sym-4-omega-1-l}
	Let $G = \Sym_4$ and $p = 2$. Let $s = 3$. Let $l = 2$.
	
	Let $U := \spn{(1,2)} \leq \Sym_4$. Then $U \in \Sub(1)$ since $|U| = 2^1$. 
	
	We choose $\tv_{\Sym_4}(U)$ such that $T := \{(3,4), (1,2,3), (1,3,4)\} \subseteq \tv_{\Sym_4}(U)$ by choice of the latter.
	
	Then $T \in \Tv_{2^2}(U)$ since $T \subseteq (\tv_{\Sym_4}(U) \setminus \{\id\})$ and $|T| = 2^2 - 1$, \cf \autoref{not-tv}.
	
	Let
	\begin{align*}
		M &:= \textstyle U \sqcup \bigsqcup_{t \in T}Ut \\
		&\phantom := U \sqcup U(3,4) \sqcup U(1,2,3) \sqcup U(1,3,4) \\
		&\phantom := \{\id, (1,2)\} \sqcup \{(3,4), (1,2)(3,4)\} \sqcup \{(1,2,3), (1,3)\} \sqcup \{(1,3,4), (1,2,3,4)\}.
	\end{align*}
	Then $M \in \Omega_1^{[0, 2]}$, \cf \autoref{definition-omega-1-l}.(2).
	
	Magma gives $\Stab_{\Sym_4}(M) = \spn{(1,2), (3,4)}$. Then $|\Stab_{\Sym_4}(M)| = 2^2 = 2^{3-1}$, \ie $$M \in \Omega_1^{1},$$ \cf \autoref{definition-omega-1-l}.(1).
\end{example}

\begin{lemma}\label{lemma-u-stab-subgroup}
	Suppose given $$M = U \sqcup \bigsqcup_{i \in [2,p^l]} Ug_i$$ for some $U \in \Sub(s-l)$ and $\{g_2, \dots, g_{p^l}\} \in \Tv_{p^l}(U)$.
	
	Then we have $$U \leq \Stab(M).$$
	
	\begin{proof}
		We show $U \stackrel!\subseteq \Stab(M)$. Suppose given $u \in U$. We show $uM \stackrel!= M$.
		
		Since $u \in U$, we have $uU = U$ and therefore \[uM
		= u\left(U \sqcup \textstyle\bigsqcup_{i \in [2,p^l]} Ug_i\right)
		= uU \sqcup \textstyle\bigsqcup_{i \in [2,p^l]} uUg_i
		= U \sqcup \textstyle\bigsqcup_{i \in [2,p^l]} Ug_i
		= M.\qedhere\]
	\end{proof}
\end{lemma}

\begin{example}
	We continue \autoref{ex-sym-4-omega-1-l}.
	
	We have $U = \spn{(1,2)}$ and $\Stab_{\Sym_4}(M) = \spn{(1,2), (3,4)}$. Indeed, we have $$U \leq \Stab_{\Sym_4}(M)$$ as predicted by \autoref{lemma-u-stab-subgroup}. 
	
	Since $|U| = 2^1$ and $|\Stab_{\Sym_4}(M)| = 2^2$, we have $$(U, \Stab_{\Sym_4}(M)) \in \Sub(1, 2).$$
\end{example}

\begin{lemma}\label{lemma-omega-0-l}
	We have $\Omega_1^{[0,l]} = \bigsqcup_{k \in [0,l]} \Omega_1^k$.
	
	\begin{proof}
		\Ad $\subseteq$. Let $$M := U \sqcup \bigsqcup_{i \in [2,p^l]} Ug_i \in \Omega_1^{[0,l]},$$ where $U \in \Sub(s-l)$ and $\{g_2, \dots, g_{p^l}\} \in \Tv_{p^l}(U)$. Then $|M| = p^{s-l} + (p^l - 1)p^{s-l} = p^s$, so $M \in \Omega$.
		
		Let $V := \Stab(M)$. \autoref{lemma-order-stab-and-orbit}.(2) yields $|V| = p^{s-k}$ for some $k \in [0,s]$. Since $U \leq V$ by \autoref{lemma-u-stab-subgroup}, $p^{s-l}$ divides $p^{s-k}$, \ie $s-l \leq s-k$, \ie $k \leq l$. Note that $M \in \Omega^k$. Since $1 \in U \subseteq M$, we have $M \in \Omega_1^k \subseteq \bigsqcup_{i \in [0,l]} \Omega_1^i$.
		
		\Ad $\supseteq$. Let $M \in \Omega_1^k$ for some $k \in [0,l]$. Let $U := \Stab(M)$. Then $|U| = |\Stab(M)| = p^{s-k}$. \autoref{lemma-order-stab-and-orbit}.(2) gives a decomposition $M = \bigsqcup_{j \in [1,p^k]} Um_j$ for some $m_1, \dots, m_{p^k} \in M$. Since $1 \in M$, the right coset $U \cdot 1$ is contained in $M$.
		
		Then we may choose $\{g_2, \dots, g_{p^k}\} \in \Tv_{p^k}(U)$ such that $$M = U \sqcup \bigsqcup_{i \in [2,p^k]} Ug_i.$$
		
		Thanks to \autoref{thm-sylow-wielandt}, we may choose $V \in \Sub(s-l)$ with $V \leq U$. Then $|U/V| = p^{s-k}/p^{s-l} = p^{l-k}$. So there exists a unique $\{g_2', \dots, g_{p^{l-k}}'\} \in \Tv_{p^{l-k}}(V)$ such that $U = V \sqcup \bigsqcup_{j \in [2, p^{l-k}]} Vg_j'$. Write $g_1 := g_1' := 1$. Then \begin{align*}
			M = \bigsqcup_{i \in [1,p^k]} Ug_i = \bigsqcup_{i \in [1,p^k]}(\textstyle\bigsqcup_{m \in [1,p^{l-k}]} Vg_m')g_i = \displaystyle\bigsqcup_{\substack{i \in [1, p^k] \\ m \in [1, p^{l-k}]}} Vg_m'g_i = V \sqcup \bigsqcup_{j \in [2, p^l]} Vh_j
		\end{align*} for some $\{h_2, \dots, h_{p^l}\} \in \Tv_{p^l}(V)$. Hence, $M \in \Omega_1^{[0,l]}$.
	\end{proof}
\end{lemma}

\begin{definition}\label{definition-data}\leavevmode
	\begin{enumerate}[listparindent=0pt]
		\item Let \begin{align*}
			\Data^{(l)} := \Data^{[s],(l)} &:= \coprod_{U \in \Sub(s-l)} \Tv_{p^l}(U) \\
			&\phantom{:}= \{(U, \{g_2, \dots, g_{p^l}\}): U \in \Sub(s-l), \{g_2, \dots, g_{p^l}\} \in \Tv_{p^l}(U)\}.
		\end{align*}
		Then $$|\Data^{(l)}| = |\Data^{[s],(l)}| = \TvBin(l)\Num(s-l);$$ \cf \autoref{remark-binomial} and \autoref{definition-sub-num}.(2).
		
		\item Let \begin{align*}
			\DataPhi^{(l)}:\quad\qquad\qquad\quad \Data^{(l)} &\to \Omega_1^{[0,l]} \\
			(U, \{g_2, \dots, g_{p^l}\}) &\mapsto U \sqcup \bigsqcup_{i \in [2,p^l]} Ug_i
		\end{align*}
		be the \textit{assembly map}. Let $$\DataRho_1^{[0,l]} := \DataRho|_{\Omega_1^{[0,l]}}^{\overline \Omega^{[0,l]}} : \Omega_1^{[0,l]} \to \overline\Omega^{[0,l]}.$$
		
		\item Let $k \in [0,l]$. Let $$\Data^{(l),k} := \Data^{[s],(l),k} := (\DataPhi^{(l)})^{-1}(\Omega_1^k) \subseteq \Data^{(l)}.$$ So $\Data^{(l)} = \bigsqcup_{i \in [0,l]} \Data^{(l),i}$, \cf \autoref{lemma-omega-0-l}.
		
		Let $$\DataPhi^{(l),k} := \DataPhi^{(l)}|_{\Data^{(l),k}}^{\Omega_1^k} : \Data^{(l),k} \to \Omega_1^k.$$ and $$\DataRho_1^k := \DataRho_1^{[0,l]}|_{\Omega_1^k}^{\overline \Omega^k} : \Omega_1^k \to \overline\Omega^k.$$
	\end{enumerate}
	The situation can be illustrated as follows.
	
	\begin{center}
		\begin{tikzcd}[column sep=1ex,row sep=5ex]
			&&\Data^{(l)}\arrow[d, "\DataPhi^{(l)}"]
			&\quad=\quad
			&\Data^{(l),0}\arrow[d, "\DataPhi^{(l),0}"]
			&\sqcup
			&\Data^{(l),1}\arrow[d, "\DataPhi^{(l),1}"]
			&\sqcup
			&\ldots
			&\sqcup
			&\Data^{(l),l}\arrow[d, "\DataPhi^{(l),l}"] \\
			
			\Omega\arrow[d, "\DataRho"]
			& \supseteq
			&\Omega_1^{[0,l]}\arrow[d, "\DataRho_1^{[0,l]}"]
			&\quad=\quad
			&\Omega_1^0\arrow[d, "\DataRho_1^{0}"]
			&\sqcup
			&\Omega_1^1\arrow[d, "\DataRho_1^{1}"]
			&\sqcup
			&\ldots
			&\sqcup
			&\Omega_1^l\arrow[d, "\DataRho_1^{l}"] \\
			
			\overline\Omega
			&\quad\supseteq\quad
			&\overline \Omega^{[0,l]}
			&=
			&\overline\Omega^0
			&\sqcup
			&\overline\Omega^1
			&\sqcup
			&\ldots
			&\sqcup
			&\overline\Omega^l
		\end{tikzcd}
	\end{center}
\end{definition}

\begin{example}
	We continue \autoref{ex-sym-4-omega-1-l}. Recall that $s = 3$.
	
	Since $U \in \Sub(3-2)$ and $T \in \Tv_{2^2}(U)$, we have $$\theta := (U, T) \in \Data^{(2)},$$ \cf \autoref{definition-data}.(1).
	
	Consider $\DataPhi^{(2)} : \Data^{(2)} \to \Omega_1^{[0, 2]}$. We have
	\begin{align*}
		\DataPhi^{(2)}(\theta) &= \DataPhi^{(2)}\bigl((U, T)\bigr) \\
		&= \DataPhi^{(2)}\bigl((U, \{(3,4), (1,2,3), (1,3,4)\})\bigr) \\
		&= \textstyle U \sqcup \bigsqcup_{g \in \{(3,4), (1,2,3), (1,3,4)\}}Ug \\
		&= U \sqcup U(3,4) \sqcup U(1,2,3) \sqcup U(1,3,4) \\
		&= M,
	\end{align*} \cf \autoref{definition-data}.(2).
	
	Since $$M = \DataPhi^{(2)}(\theta)\ \ \stackrel{\text{Ex. \ref{ex-sym-4-omega-1-l}}}{\in}\ \Omega_1^{1}\ \subseteq\ \Omega_1^{0} \sqcup \Omega_1^{1} \sqcup \Omega_1^{2}\ \stackrel{\text{L. \ref{lemma-omega-0-l}}}{=}\ \Omega_1^{[0,2]},$$ we have $$(U, T) \in (\DataPhi^{(2)})^{-1}(\{M\}) = \Data^{(2),1},$$ and $$(U, T)\ \stackrel{\DataPhi^{(2),1}}\longmapsto\ M,$$ \cf \autoref{definition-data}.(3).
\end{example}

%% file: text/BA-Some_basic_properties.tex
Let $l \in [0, s]$.

\subsection{A characterization}
\begin{lemma}\label{lemma-phi-l-surjective}\leavevmode
	\begin{enumerate}[listparindent=0pt]
		\item The map $\DataPhi^{(l)}: \Data^{(l)} \to \Omega_1^{[0,l]}$ is surjective.
		\item The map $\DataPhi^{(l),k}: \Data^{(l),k} \to \Omega_1^k$ is surjective for $k \in [0, l]$.
	\end{enumerate}
	
	\begin{proof}
		\Ad (1). \Cf \autoref{definition-omega-1-l}.(2) and \autoref{definition-data}.(1).
		
		\Ad (2). \Cf (1) and \autoref{definition-data}.(2, 3).
	\end{proof}
\end{lemma}

\begin{lemma}\label{lemma-data-equivalences}
	Let $(U, \{g_2, \dots, g_{p^l}\}) \in \Data^{(l)}$. Let $$M := U \sqcup \bigsqcup_{i \in [2,p^l]} Ug_i = \DataPhi^{(l)}(U, \{g_2, \dots, g_{p^l}\}) \in \Omega_1^{[0, l]}.$$ The following assertions are equivalent.
	\begin{enumerate}[listparindent=0pt]
		\item We have $(U, \{g_2, \dots, g_{p^l}\}) \in \Data^{(l),l}$.
		\item We have $M \in \Omega_1^l$.
		\item We have $U = \Stab(M)$.
	\end{enumerate}
	
	\begin{proof}
		\Ad (1) $\Leftrightarrow$ (2). By \autoref{definition-data}.(2), we have $(U, \{g_2, \dots, g_{p^l}\}) \in \Data^{(l),l}$ if and only if $\DataPhi^{(l)}(U, \{g_2, \dots, g_{p^l}\}) \in \Omega_1^l$.
		
		\Ad (2) $\Rightarrow$ (3). We have $M = U \sqcup \bigsqcup_{i \in [2,p]} Ug_i \in \Omega_1^l$.
		\autoref{lemma-u-stab-subgroup} yields $U \leq \Stab(M)$. We have $|U| = p^{s-l}$ by \autoref{definition-data}.(1). We have $|\Stab(M)| = p^{s-l}$ by \autoref{definition-omega}.(1). This yields $U = \Stab(M)$.
		
		\Ad (3) $\Rightarrow$ (2). If $U = \Stab(M)$, then $|\Stab(M)| = |U| = p^{s-l}$; \cf \autoref{definition-data}.(1). So $M \in \Omega^l$; \cf \autoref{definition-omega}.(1). Since $1 \in U \subseteq M$, we get $M \in \Omega_1^l$.			
	\end{proof}
\end{lemma}

\begin{example}
	We continue \autoref{ex-sym-4-omega-1-l}. Recall that $s = 3$.
	
	We have $$U = \spn{(1,2)} < \spn{(1,2), (3,4)} = \Stab_{\Sym_4}(M),$$ and therefore $M \notin \Omega_1^{2}$, \cf \autoref{lemma-data-equivalences}.($\lnot 3 \Rightarrow \lnot 2$).
	
	Actually, we have already seen that $M \in \Omega_1^{1}$; cf. \autoref{ex-sym-4-omega-1-l}.
	
	Now let $l = 1$ and consider $V := \Stab_{\Sym_4}(M)$. Suppose that $(1,2,3) \in \tv_{\Sym_4}(V)$ by choice of the latter. Then $V \in \Sub(3-1)$ and $\{(1,2,3)\} \in \Tv_2(V)$, and we get $$M = V \sqcup V(1,2,3) = \DataPhi^{(1)}\big((V, \{(1,2,3)\})\big).$$ The equivalent assertions of \autoref{lemma-data-equivalences} all hold: We have $(V, \{(1,2,3)\}) \in \Data^{(1),1}$, $M \in \Omega_1^1$ and $V = \Stab_{\Sym_4}(M)$. In particular, $M = \DataPhi^{(1),1}\big((V, \{(1,2,3)\})\big)$.
\end{example}

\subsection{Three bijections}
\begin{lemma}\label{lemma-bijection-sub-sml}
	We have the bijective map \begin{align*}
		\Data^{(l),0} &\to \Sub(s-l, s) \\
		(U, \{g_2, \dots, g_{p^l}\}) &\mapsto (U, \textstyle U \sqcup \bigsqcup_{i \in [2, p^l]} Ug_i).
	\end{align*} In particular, we have $$|\Data^{(l),0}| = \Num(s-l, s),$$ \cf \autoref{definition-sub-num}\upshape.(2).
	\begin{proof}
		\emph{Well-defined}. Given $(U, \{g_2, \dots, g_{p^l}\}) \in \Data^{(l),0}$, we have $$\DataPhi^{(l),0}\bigl((U, \{g_2, \dots, g_{p^l}\})\bigr) = U \sqcup \textstyle\bigsqcup_{i \in [2, p^l]} Ug_i =: M \in \Omega_1^0.$$ Let $V := \Stab(M)$. We have $|V| = p^s$; \cf \autoref{definition-omega}.(1). We have $U \leq V$; \cf \autoref{lemma-u-stab-subgroup}. Since $|[M]| = p^{t-s}n$, we get $M = Vm$ for some $m \in M$; \cf \autoref{lemma-order-stab-and-orbit}.(2). But $1 \in M = Vm$, so $V1 = Vm$ and therefore $M = V \leq G$. Hence $(U, V) \in \Sub(s-l, s)$.
		
		\emph{Surjective}. Suppose given $(U, V) \in \Sub(s-l, s)$. We have $|V|/|U| = p^s/p^{s-l} = p^l$, so we may choose $\{g_2, \dots, g_{p^l}\} \in \Tv_{p^l}(U)$ with $U \sqcup \bigsqcup_{i \in [2, p^l]} Ug_i = V$. Then $V \in \Omega_1^0$, \cf \autoref{remark-stab-of-subgroup}. Since $$V = \DataPhi^{(l)}\bigl((U, \{g_2, \dots, g_{p^l}\})\bigr) \in \Omega_1^0,$$ we have $(U, \{g_2, \dots, g_{p^l}\}) \in \Data^{(l),0}$.
		
		\emph{Injective}. Given $(U, \{g_2, \dots, g_{p^l}\}), (\hat U, \{\hat g_2, \dots, \hat g_{p^l}\}) \in \Data^{(l),0}$ with $$(U, \textstyle U \sqcup \bigsqcup_{i \in [2, p^l]} Ug_i) = (\hat U, \textstyle \hat U \sqcup \bigsqcup_{i \in [2, p^l]} \hat U\hat g_i),$$ it follows first that $U = \hat U$.
		
		Since $\{g_2, \dots, g_{p^l}\}, \{\hat g_2, \dots, \hat g_{p^l}\} \in \Tv_{p^l}(U)$ and $U \sqcup \bigsqcup_{i \in [2, p^l]} Ug_i = U \sqcup \bigsqcup_{i \in [2, p^l]} U\hat g_i$, we conclude that $$\{g_2, \dots, g_{p^l}\} = \{\hat g_2, \dots, \hat g_{p^l}\},$$ \cf \autoref{remark-disjoint-union-index-sets}.
	\end{proof}
\end{lemma}

\begin{lemma}\label{lemma-bijection-omega-1-l}
	The map $$\DataPhi^{(l),l} : \Data^{(l),l} \to \Omega_1^l,\ (U, \{g_2, \dots, g_{p^l}\}) \mapsto U \textstyle\sqcup \bigsqcup_{i \in [2, p^l]} Ug_i$$ is bijective. Moreover, $$U = \Stab(\textstyle U \sqcup \bigsqcup_{i \in [2, p^l]} Ug_i)$$ for $(U, \{g_2, \dots, g_{p^l}\}) \in \Data^{(l),l}$, \cf \autoref{lemma-data-equivalences}.$(2\Rightarrow3)$.

	\begin{proof}
		\emph{Surjective}. \Cf \autoref{lemma-phi-l-surjective}.(2).
		
		\emph{Injective}. Let $(U, \{g_2, \dots, g_{p^l}\}), (\hat U, \{\hat g_2, \dots, \hat g_{p^l}\}) \in \Data^{(l),l}$ with $$M := \DataPhi^{(l),l}(U, \{g_2, \dots, g_{p^l}\}) = \DataPhi^{(l),l}(\hat U, \{\hat g_2, \dots, \hat g_{p^l}\}),$$ \ie $$U \sqcup \bigsqcup_{i \in [2, p^l]} Ug_i = \hat U \sqcup \bigsqcup_{i \in [2, p^l]} \hat U\hat g_i.$$ Then $U = \Stab(M) = \hat U$, \cf \autoref{lemma-data-equivalences}.($2 \Rightarrow 3$). Since $\Tv_{p^l}(U) = \Tv_{p^l}(\hat U)$ and since $\{Ug_2, \dots, Ug_{p^l}\} = \{U\hat g_2, \dots, U\hat g_{p^l}\}$, it follows that $$\{g_2, \dots, g_{p^l}\} = \{\hat g_2, \dots, \hat g_{p^l}\},$$ \cf \autoref{remark-disjoint-union-index-sets}. So \[(U, \{g_2, \dots, g_{p^l}\}) = (\hat U, \{\hat g_2, \dots, \hat g_{p^l}\}).\qedhere\]
	\end{proof}
\end{lemma}

\begin{lemma}\label{lemma-bijection-sub-sml-stab}
	Let $k \in [0, l]$. Let $M \in \Omega_1^k$. We have the bijective map
	\begin{align*}
		\kappa:\qquad (\DataPhi^{(l),k})^{-1}(\{M\}) &\to \Sub(s-l, \Stab(M)) \\
		(V, \{h_2, \dots, h_{p^l}\}) &\mapsto V
	\end{align*} In particular, we have $$|(\DataPhi^{(l),k})^{-1}(\{M\})| =  |\Sub(s-l, \Stab(M))|.$$
	
	\begin{proof}
		There exists a unique $(U, \{g_2, \dots, g_{p^k}\}) \in \Data^{(k),k}$ such that $$M = U \sqcup \bigsqcup_{i \in [2, p^k]} Ug_i,$$ where $U = \Stab(M)$, \cf \autoref{lemma-bijection-omega-1-l}.
		
		\emph{Well-defined}. Suppose given $(V, \{h_2, \dots, h_{p^l}\}) \in (\DataPhi^{(l),k})^{-1}(\{M\}) \subseteq \Data^{(l),k}$. Then $$V \sqcup \textstyle\bigsqcup_{j \in [2, p^l]}Vh_j = M = U \sqcup \bigsqcup_{i \in [2, p^k]} Ug_i,$$ and, by \autoref{lemma-u-stab-subgroup}, $V \leq \Stab(M) = U$, so $V \in \Sub(s-l, U)$.
		
		\emph{Surjective}. Suppose given $V \in \Sub(s-l)$ with $V \leq U$.
		
		There exists a unique $\{g_2', \dots, g_{p^{l-k}}'\} \in \Tv_{p^{l-k}}(V)$ such that $U = V \sqcup \bigsqcup_{m \in [2, p^{l-k}]} Vg_m'$.
		
		Write $g_1 := g_1' := 1$. We get
		\begin{align*}
			M = \bigsqcup_{i \in [1, p^k]} Ug_i 
			= \bigsqcup_{i \in [1, p^k]}(\textstyle\bigsqcup_{m \in [1, p^{l-k}]} Vg_m')g_i
			= \displaystyle\bigsqcup_{\substack{i \in [1, p^k] \\ m \in [1, p^{l-k}]}} Vg_m'g_i 
			= V \sqcup \bigsqcup_{j \in [2, p^l]} Vh_j
		\end{align*} for some $\{h_2, \dots, h_{p^l}\} \in \Tv_{p^l}(V)$.
		
		Since $M \in \Omega_1^k$, we have $(V, \{h_2, \dots, h_{p^l}\}) \in \Data^{(l),k}$.
		
		Moreover, since $M = V \sqcup \bigsqcup_{j \in [2, p^l]} Vh_j = \DataPhi^{(l),k}\bigl((V, \{h_2, \dots, h_{p^l}\})\bigr)$, we have $$(V, \{h_2, \dots, h_{p^l}\}) \in (\DataPhi^{(l),k})^{-1}(\{M\}).$$
		
		Hence $$V = \kappa\big((V, \{h_2, \dots, h_{p^l}\})\big) \in \kappa\big((\DataPhi^{(l),k})^{-1}(\{M\})\big).$$
		\begin{note}
			Note that we had a similar calculation in \autoref{lemma-omega-0-l}, where $\bigsqcup_{k \in [0,l]} \Omega_1^k \stackrel !\subseteq \Omega_1^{[0,l]}$ was shown.
		\end{note}
		
		\emph{Injective}. Given $(V, \{h_2, \dots, h_{p^l}\}), (\hat V, \{\hat h_2, \dots, \hat h_{p^l}\}) \in (\DataPhi^{(l),k})^{-1}(\{M\})$ which map to the same element $V = \hat V$ under $\kappa$, we have $$V \sqcup \textstyle\bigsqcup_{j \in [2, p^l]}Vh_j = M = V \sqcup \textstyle\bigsqcup_{j \in [2, p^l]}V\hat h_j.$$
		Therefore, $\{h_2, \dots, h_{p^l}\} = \{\hat h_2, \dots, \hat h_{p^l}\}$, \cf \autoref{remark-disjoint-union-index-sets}.
	\end{proof}
\end{lemma}

\begin{example}
	Let $G = \Sym_4$ and $p = 2$. Let $s = 3$.
	
	\begin{enumerate}[listparindent=0pt]
		\item\label{ex-m-finden-teil-1}
		Let $l = 1$. We want to find an element $$M \in \Omega_1^1 = \Omega_1^{[3],1} = \{M \subseteq \Sym_4 : |M| = 2^3,\, |\Stab_{\Sym_4}(M)| = 2^{3-1}\}$$ by making use of \autoref{lemma-bijection-omega-1-l}.
		
		Recall that $$\Data^{(1)} = \Data^{[3],(1)} = \{(U, \{g_2\}) : U \in \Sub(3-1),\, \{g_2\} \in \Tv_{2^1}(U)\},$$ \cf \autoref{definition-data}.(1).
		
		Recall that we have $\Data^{(1)} = \Data^{(1),0} \sqcup \Data^{(1),1}$, \cf \autoref{definition-data}.(3).
		
		Suppose given $(U, \{g_2\}) \in \Data^{(1)}$, recall that we have $(U, \{g_2\}) \in \Data^{(1),1}$ if and only if $\DataPhi^{(1)}\bigl((U, \{g_2\})\bigr) = U \sqcup Ug_2 \nleq \Sym_4$, \cf \autoref{lemma-bijection-sub-sml} and \autoref{lemma-wielandt}.(2).
		
		Let $U := \{\id, (1,2), (3,4), (1,2)(3,4)\} = \spn{(1,2), (3,4)} \in \Sub(3-1)$.
		
		We suppose that $(2,3) \in \tv_{\Sym_4}(U)$ by choice of the latter.
		
		So $(U, \{(2,3)\}) \in \Data^{(1)}$.
		
		Let
		\begin{align*}
			M &:= \DataPhi^{(1)}\bigl((U, \{g_2\})\bigr) = U \sqcup U(2,3) \\
			&\phantom := \{\id, (1,2), (3,4), (1,2)(3,4), (2,3), (1,3,2), (2,3,4), (1,3,4,2)\}.
		\end{align*}
		Then $M \nleq G$ since $(1,2,3) \in M$ and $|\spn{(1,2,3)}| = 3$ does not divide $|M| = 8$.
		
		Thus, $(U, \{(2,3)\}) \in \Data^{(1),1}$, \ie $M \in \Omega_1^1$.
		
		Furthermore, \autoref{lemma-bijection-omega-1-l} gives $\Stab_{\Sym_4}(M) = U$.
		
		\item
		Now let $l = 2$. Let $M \in \Omega_1^1$ be as in \ref{ex-m-finden-teil-1}.
		
		We want to consider the bijection
		\begin{align*}
			\kappa:\qquad (\DataPhi^{(2),1})^{-1}(\{M\}) &\to \Sub(3-2, U) \\
			(V, \{h_2, h_3, h_4\}) &\mapsto V
		\end{align*}
		as in \autoref{lemma-bijection-sub-sml-stab}, where $k = 1$.
		
		Let $V_1 := \spn{(1,2)}$, $V_2 := \spn{(3,4)}$ and $V_3 := \spn{(1,2)(3,4)}$.
		
		Then $V_1, V_2, V_3 \in \Sub(3-2)$, and we have $\Sub(3-2, U) = \{V_1, V_2, V_3\}$.
		
		We have $M = V_1 \sqcup V_1(3,4) \sqcup V_1(2,3) \sqcup V_1(2,3,4)$.
		
		Supposing $\{(3,4), (2,3), (2,3,4)\} \subseteq \tv_{\Sym_4}(V_1)$ by choice of the latter, we have $$(V_1, \{(3,4), (2,3), (2,3,4)\}) \stackrel \kappa\mapsto V_1.$$
		
		We have $M = V_2 \sqcup V_2(1,2) \sqcup V_2(2,3) \sqcup V_2(1,3,2)$.
		
		Supposing $\{(1,2), (2,3), (1,3,2)\} \subseteq \tv_{\Sym_4}(V_2)$ by choice of the latter, we have $$(V_2, \{(1,2), (2,3), (1,3,2)\}) \stackrel \kappa\mapsto V_2.$$
		
		We have $M = V_3 \sqcup V_3(1,2) \sqcup V_3(2,3) \sqcup V_3(1,3,2)$.
		
		Supposing $\{(1,2), (2,3), (1,3,2)\} \subseteq \tv_{\Sym_4}(V_3)$ by choice of the latter, we have $$(V_3, \{(1,2), (2,3), (1,3,2)\}) \stackrel \kappa\mapsto V_3.$$
%	Let $N := \{\id, (1,2), (3,4), (1,2)(3,4)\} = \spn{(1,2), (3,4)} \in \Omega \cap \Sub(2)$.
%	
%	Let $V_1 := \spn{(1,2)}$, $V_2 := \spn{(3,4)}$ and $V_3 := \spn{(1,2)(3,4)}$. Then $V_1, V_2, V_3 \in \Sub(1)$.
%	
%	Suppose that $(3,4) \in \tv_{\Sym_4}(V_1)$ and $(1,2) \in \tv_{\Sym_4}(V_2) \cap \tv_{\Sym_4}(V_3)$ by choice of the latter.
%	
%	We have $\Stab_{\Sym_4}(N) = N$ and $N \in \Omega_1^{0}$, \cf \autoref{remark-stab-of-subgroup}.
%	
%	We have $\Stab_{\Sym_4}(N) = N$ again by \autoref{lemma-bijection-omega-1-l}.
%	
%	Furthermore, we have $$N
%	= V_1 \sqcup V_1(3,4)
%	= V_2 \sqcup V_2(1,2)
%	= V_3 \sqcup V_3(1,2),$$ \ie $$N
%	= \DataPhi^{(1),0}\bigl((V_1, \{(3,4)\})\bigr)
%	= \DataPhi^{(1),0}\bigl((V_2, \{(1,2)\})\bigr)
%	= \DataPhi^{(1),0}\bigl((V_3, \{(1,2)\})\bigr).$$
%	By \autoref{ex-sym-4-sub-num}, we know that $$\Sub(2, N) \stackrel{\text{D. \ref{definition-sub-num}.(3)}}{=} \{(V, N) : V \leq G, V \subseteq N\} = \{(V_1, N), (V_2, N), (V_3, N)\}.$$ Hence $$(\DataPhi^{(1), 0})^{-1}(\{M\}) = \{(V_1, \{(3,4)\}), (V_2, \{(1,2)\}), (V_3, \{(1,2)\})\},$$ \cf \autoref{lemma-bijection-sub-sml-stab}.
	\end{enumerate}
\end{example}

%% file: text/BA-a_l.tex
\begin{lemma}\label{lemma-d-k-l-preparation}\leavevmode
	Let $l \in [0, s]$. Let $k \in [0, l]$. Let $\hat I \subseteq [l, s]$. We have $$\sum_{(U, \{g_2, \dots, g_{p^l}\}) \in \Data^{(l),k}} |\Sub(s-\hat I, U)| = \sum_{\substack{I \subseteq [0, k] \\ k \in I}} (-1)^{|I| + 1} \cdot \TvBin(\min I) \Num(s - (I \cup \{l\} \cup \hat I)).$$ \Cf Definitions {\upshape \ref{definition-sub-num}.(2, 3, 4)} and {\upshape \ref{definition-binomial}}.
	
	\begin{proof}
		Let $X := \{(k, l) \in \Z_{\geq 0} \times \Z_{\geq 0}: k \leq l \leq s\}$. For $(k, l), (k', l') \in X$, we define: $$(k, l) \preccurlyeq (k', l') \ :\Leftrightarrow\ k \leq k' \text{ and } l \leq l'.$$
		Then $(X, \preccurlyeq)$ is a partially ordered set. By \autoref{cor-finite-posets-noetherian}, $(X, \preccurlyeq)$ is noetherian.
		
		Therefore, we can use noetherian induction over $(X, \preccurlyeq)$ to prove the lemma, \cf \autoref{lemma-noetherian-induction}.
%		For $k = 0$, we get
%		\[\def\arraystretch{1.5}
%			\begin{array}{lccl}
%				&\mystackrel{}{\phantom =}
%				&\displaystyle\sum_{(U, \{g_2, \dots, g_{p^l}\}) \in \Data^{(l),0}} &|\Sub(s-\hat I, U)| \\
%				&\mystackrel{L. \ref{lemma-bijection-sub-sml}}{=}
%				&\displaystyle\sum_{(U, V) \in \Sub(s-l, s)} &|\Sub(s-\hat I, U)| \\
%				&\mystackrel{}{=}&&\Num(s-(\{0\} \cup \{l\} \cup \hat I)) \\
%				&\mystackrel{}{=}
%				&\displaystyle\sum_{\substack{I \subseteq [0, 0] \\ 0 \in I}}
%				&(-1)^{|I| + 1} \cdot \TvBin(\min I) \Num(s - (I \cup \{l\} \cup \hat I)).
%			\end{array}
%		\]

		Let $(k, l) \in X$. Suppose that the claim is true for every $(i, j) \in X$ with $(i, j) \prec (k, l)$. We obtain
		
		\allowdisplaybreaks
		\[\def\arraystretch{1.5}
			\begin{array}{ccccl}
				\mystackrel{}{\phantom =}
				&&&\displaystyle\sum_{(U, \{g_2, \dots, g_{p^l}\}) \in \Data^{(l),k}} &|\Sub(s-\hat I, U)| \\
				\mystackrel{D. \ref{definition-data}.(3)}{=}
				&&\displaystyle\sum_{M \in \Omega_1^k}
				&\displaystyle\sum_{(U, \{g_2, \dots, g_{p^l}\}) \in (\DataPhi^{(l),k})^{-1}(\{M\})}
				&|\Sub(s-\hat I, U)| \\
				\mystackrel{L. \ref{lemma-bijection-sub-sml-stab}}{=}
				&&\displaystyle\sum_{M \in \Omega_1^k}
				&\displaystyle\sum_{U \in \Sub(s-l, \Stab(M))}
				&|\Sub(s-\hat I, U)| \\
				\mystackrel{L. \ref{lemma-bijection-omega-1-l}}{=}
				&&\displaystyle\sum_{(V, \{h_2, \dots, h_{p^k}\}) \in \Data^{(k),k}}
				&\displaystyle\sum_{U \in \Sub(s-l, V)}
				&|\Sub(s-\hat I, U)|
				% ACHTUNG: Hier wurde der Zeilenumbruch hartkodiert!
			\end{array}
		\]
		\[\def\arraystretch{1.5}
			\begin{array}{ccccl}
				\mystackrel{}{=}
				&&&\displaystyle\sum_{(V, \{h_2, \dots, h_{p^k}\}) \in \Data^{(k),k}}
				&|\Sub(s-(\{l\} \cup \hat I), V)| \\
				\mystackrel{D. \ref{definition-data}.(3)}{=}
				&&&\displaystyle\sum_{(V, \{h_2, \dots, h_{p^k}\}) \in \Data^{(k)}}
				&|\Sub(s-(\{l\} \cup \hat I), V)| \\
				\mystackrel{}{\phantom{=}}
				&-&\displaystyle\sum_{j \in [0, k-1]}
				&\displaystyle\sum_{(V, \{h_2, \dots, h_{p^k}\}) \in \Data^{(k),j}}
				&|\Sub(s-(\{l\} \cup \hat I), V)| \\
				\mystackrel{D. \ref{definition-data}.(1), R. \ref{remark-binomial}}{=}
				&&&&\TvBin(k)\Num(s-(\{k, l\} \cup \hat I)) \\
				\mystackrel{}{\phantom{=}}
				&-&\displaystyle\sum_{j \in [0, k-1]}
				&\displaystyle\sum_{(V, \{h_2, \dots, h_{p^k}\}) \in \Data^{(k),j}}
				&|\Sub(s-(\{l\} \cup \hat I), V)|.
			\end{array}
		\]

		For $j \in [0, k-1]$, we have $(j, k) \prec (k, l)$ and $\{l\} \cup \hat I \subseteq [k, s]$, and so $$\sum_{(V, \{h_2, \dots, h_{p^k}\}) \in \Data^{(k),j}} |\Sub(s-(\{l\} \cup \hat I), V)| = \sum_{\substack{I \subseteq [0, j] \\ j \in I}} (-1)^{|I| + 1} \TvBin(\min I)\Num(s - (I \sqcup (\{k, l\} \cup \hat I)))$$ by induction hypothesis. Hence
		\[\def\arraystretch{1.5}
			\begin{array}{lcccl}
				&\mystackrel{}{\phantom{=}}
				&\displaystyle\sum_{j \in [0, k-1]}
				&\displaystyle\sum_{(V, \{h_2, \dots, h_{p^k}\}) \in \Data^{(k),j}}
				&|\Sub(s-(\{l\} \cup \hat I), V)| \\
				&\mystackrel{\text{IH}}{=}
				&\displaystyle\sum_{j \in [0, k-1]}
				&\displaystyle\sum_{\substack{I \subseteq [0, j] \\ j \in I}}
				&(-1)^{|I| + 1} \TvBin(\min I)\Num(s - (I \sqcup (\{k, l\} \cup \hat I))) \\
				&\mystackrel{}{=}
				&&\displaystyle\sum_{\substack{I \subseteq [0, k-1] \\ I \neq \emptyset}}
				&(-1)^{|I| + 1} \TvBin(\min I)\Num(s - (I \sqcup (\{k, l\} \cup \hat I))) \\
				&\mystackrel{$I' = I \cup \{k\}$}{=}
				&&\displaystyle\sum_{\substack{I' \subseteq [0, k] \\ I' \neq \{k\},\ k \in I'}}
				&(-1)^{|I'|} \TvBin(\min I')\Num(s - (I' \cup \{l\} \cup \hat I)). \\
			\end{array}
		\]
		Therefore,
		\[\def\arraystretch{1.5}
			\begin{array}{lcccl}
				\mystackrel{}{\phantom{=}}
				&&&&\TvBin(k)\Num(s-(\{k, l\} \cup \hat I)) \\
				\mystackrel{}{\phantom{=}}
				&-
				&\displaystyle\sum_{j \in [0, k-1]}
				&\displaystyle\sum_{(V, \{h_2, \dots, h_{p^k}\}) \in \Data^{(k),j}}
				&|\Sub(s-(\{l\} \cup \hat I), V)| \\
				\mystackrel{}{=}
				&&&&\TvBin(k)\Num(s-(\{k\} \cup \{l\} \cup \hat I)) \\
				\mystackrel{}{\phantom{=}}
				&-
				&&\displaystyle\sum_{\substack{I' \subseteq [0, k] \\ I' \neq \{k\},\ k \in I'}}
				&(-1)^{|I'|} \TvBin(\min I')\Num(s - (I' \cup \{l\} \cup \hat I)) \\
				\mystackrel{}{=}
				&
				&&\displaystyle\sum_{\substack{I \subseteq [0, k] \\ k \in I}}
				&(-1)^{|I|+1} \TvBin(\min I)\Num(s - (I \cup \{l\} \cup \hat I)), \\
			\end{array}
		\]
		which completes the induction.
	\end{proof}
\end{lemma}

\begin{corollary}\label{cor-cardinality-d-k-l}
	Let $l \in [0, s]$. Let $k \in [0, l]$. We have $$|\Data^{(l),k}| = \sum_{\substack{I \subseteq [0, k] \\ k \in I}}
	(-1)^{|I| + 1} \cdot \TvBin(\min I) \Num(s - (I \cup \{l\}))$$
	
	\begin{proof}
		We have
		\begin{alignat*}{4}
			&&\qquad\qquad&\qquad\qquad|\Data^{(l),k}| \\
			&\mystackrel{}{=}
			&\smashoperator[rl]{\sum_{(U, \{g_2, \dots, g_{p^l}\}) \in \Data^{(l),k}}}
			&\qquad\qquad1 \\
			&\mystackrel{$|U| = p^{s-l}$}{=}
			&\smashoperator[rl]{\sum_{(U, \{g_2, \dots, g_{p^l}\}) \in \Data^{(l),k}}}
			&\qquad\qquad|\Sub(s-\{l\}, U)| \\
			&\mystackrel{L. \ref{lemma-d-k-l-preparation}}{=}
			&\smashoperator[rl]{\sum_{\substack{I \subseteq [0, k] \\ k \in I}}}
			&\qquad\qquad(-1)^{|I| + 1} \cdot \TvBin(\min I) \Num(s - (I \cup \{l\} \cup \{l\})) \\
			&\mystackrel{}{=}
			&\smashoperator[rl]{\sum_{\substack{I \subseteq [0, k] \\ k \in I}}}
			&\qquad\qquad(-1)^{|I| + 1} \cdot \TvBin(\min I) \Num(s - (I \cup \{l\})).
			\qedhere
		\end{alignat*}
	\end{proof}
\end{corollary}

\begin{lemma}\label{lemma-fibres-phi-bar}
	Let $l \in [0, s]$. Let $M \in \Omega_1^l$.
	
	Let $$\overline \DataPhi^{(l),l} := \DataRho_1^{l} \circ \DataPhi^{(l),l} : \Data^{(l),l} \to \overline \Omega^l.$$ Note that
	\begin{align*}
		(U, \{g_2, \dots, g_{p^l}\})\
		&\stackrel{\DataPhi^{(l),l}}{\longmapsto}\
		\textstyle U \sqcup \bigsqcup_{i \in [2, p^l]} Ug_i \\
		&\stackrel{\DataRho_1^l}{\longmapsto}\
		\textstyle [U \sqcup \bigsqcup_{i \in [2, p^l]} Ug_i] = \overline\DataPhi^{(l),l}\bigl((U, \{g_2, \dots, g_{p^l}\})\bigr)
	\end{align*}
	for $(U, \{g_2, \dots, g_{p^l}\}) \in \Data^{(l)}$.
	
	The map $\overline \DataPhi^{(l),l}$ is surjective; \cf \autoref{lemma-phi-l-surjective}{\upshape.(2)}.
	\begin{enumerate}[listparindent=0pt]
		\item\label{lemma-fibres-phi-bar-1} We have $|(\overline \DataPhi^{(l),l})^{-1}(\{[M]\})| = p^l$.
		\item We have $a_l = |\overline\Omega^l| = \tfrac 1{p^l}|\Data^{(l),l}|$, \cf \autoref{definition-omega}\upshape.(3).
	\end{enumerate}
	
	\begin{proof}
		\Ad (1). By \autoref{lemma-bijection-omega-1-l}, there exists a unique $(U, \{g_2, \dots, g_{p^l}\}) \in \Data^{(l),l}$ such that $$M = U \sqcup \bigsqcup_{i \in [2,p^{l}]} Ug_i,$$
		where $U = \Stab(M)$.
		
		In the following, let us write $g_1 := 1 \in G$ and $r := p^l \in \Z_{\geq 1}$. 
		
		We define the mapping\begin{align*}
			\lambda: [1,r] &\to (\overline \DataPhi^{(l),l})^{-1}(\{[M]\}) \\
			k &\mapsto ({}^{g_k^{-1}}U, \{g_{k,2}, \dots, g_{k,r}\})
		\end{align*}
		where $\{g_{k,2}, \dots, g_{k,r}\} \in \Tv_{r}({}^{g_k^{-1}}U)$ is such that $$\{{}^{g_k^{-1}}Ug_{k,2}, \dots, {}^{g_k^{-1}}Ug_{k,r}\} = \{{}^{g_k^{-1}}Ug_k^{-1}g_j : j \in [1,r] \setminus \{k\}\}$$ holds.
		
		We \emph{claim} that $\lambda$ is a well-defined, bijective map.
		
		\emph{Well-defined}. Let $k \in [1,r]$.
		\begin{enumerate}[label=(\roman*),listparindent=0pt]
			\item We recall that $(\overline \DataPhi^{(l),l})^{-1}(\{[M]\}) \subseteq \Data^{(l)} = \coprod_{U \in \Sub(s-l)} \Tv_{p^l}(U)$, \cf \autoref{definition-data}.(1).
			
			\begin{enumerate}[listparindent=0pt]
				\item We have ${}^{g_k^{-1}}U \in \Sub(s-l)$ since $|{}^{g_k^{-1}}U| = |U| = p^{s-l}$.
				
				\item Let $j, j' \in [1,r] \setminus \{k\}$ with $j \neq j'$. We have to show ${}^{g_k^{-1}}Ug_k^{-1}g_j \cap {}^{g_k^{-1}}Ug_k^{-1}g_{j'} \stackrel != \emptyset$.
				
				We obtain $g_k^{-1}Ug_j = {}^{g_k^{-1}}Ug_k^{-1}g_j$ and $g_k^{-1}Ug_{j'} = {}^{g_k^{-1}}Ug_k^{-1}g_{j'}$. Since $Ug_j \cap Ug_{j'} = \emptyset$, we have $g_k^{-1}Ug_j \cap g_k^{-1}Ug_{j'} = \emptyset$.
				
				\item Let $j \in [1,r] \setminus \{k\}$. We have to show ${}^{g_k^{-1}}Ug_k^{-1}g_j \cap {}^{g_k^{-1}}U \stackrel != \emptyset$.
				
				We obtain $g_k^{-1}Ug_j = {}^{g_k^{-1}}Ug_k^{-1}g_j$ and $g_k^{-1}Ug_k = {}^{g_k^{-1}}U$. Since $j \neq k$, we have $Ug_j \cap Ug_k = \emptyset$ and therefore $g_k^{-1}Ug_j \cap g_k^{-1}Ug_k = \emptyset$.
			\end{enumerate}
			
			By (a, b, c), we see that the elements of $\{{}^{g_k^{-1}}U\} \sqcup \{{}^{g_k^{-1}}Ug_k^{-1}g_j : j \in [1,r] \setminus \{k\}\}$ are pairwise distinct right-cosets of ${}^{g_k^{-1}}U$.
			
			Therefore, there is a unique $\{g_{k,2}, \dots, g_{k,r}\} \in \Tv_{r}({}^{g_k^{-1}}U)$ satisfying $$\{{}^{g_k^{-1}}Ug_{k,2}, \dots, {}^{g_k^{-1}}Ug_{k,r}\} = \{{}^{g_k^{-1}}Ug_k^{-1}g_j : j \in [1,r] \setminus \{k\}\}.$$			
			All in all, we have shown that $({}^{g_k^{-1}}U, \{g_{k,2}, \dots, g_{k,r}\}) \in \Data^{(l)}$.
			
			\item We have to show that $({}^{g_k^{-1}}U, \{g_{k,2}, \dots, g_{k,r}\}) \in (\overline \DataPhi^{(l),l})^{-1}(\{[M]\})$,
			\ie $$\overline \DataPhi^{(l),l}\bigl(({}^{g_k^{-1}}U, \{g_{k,2}, \dots, g_{k,r}\})\bigr) = [M],$$
			\ie $$[{}^{g_k^{-1}}U \sqcup \textstyle\bigsqcup_{i \in [2,r]} {}^{g_k^{-1}}Ug_{k,i}] = [M].$$
			
			We obtain \begin{align*}
				{}^{g_k^{-1}}U \sqcup \textstyle\bigsqcup_{i \in [2,r]} {}^{g_k^{-1}}Ug_{k,i}
				&= {}^{g_k^{-1}}U \sqcup \textstyle\bigsqcup_{j \in [1,r] \setminus \{k\}} {}^{g_k^{-1}}Ug_k^{-1}g_j \\
				&= g_k^{-1}(Ug_k \sqcup \textstyle\bigsqcup_{j \in [1,r] \setminus \{k\}} Ug_j) \\
				&= g_k^{-1}(U \sqcup \textstyle\bigsqcup_{i \in [2,r]} Ug_i) \\
				&= g_k^{-1}M.
			\end{align*}
			Hence $$[{}^{g_k^{-1}}U \sqcup \textstyle\bigsqcup_{i \in [2,r]} {}^{g_k^{-1}}Ug_{k,i}] = [g_k^{-1}M] = [M].$$
		\end{enumerate}
		
		\emph{Surjective}. Let $(\hat U, \{\hat g_2, \dots, \hat g_r\}) \in (\overline \DataPhi^{(l),l})^{-1}(\{[M]\}) \subseteq \Data^{(l)}$. % Hier wurde \{\} hinzugefügt.
		Then we have $$\overline \DataPhi^{(l),l}\bigl((\hat U, \{\hat g_2, \dots, \hat g_r\})\bigr) = [M] = \overline \DataPhi^{(l),l}\bigl((U, \{g_2, \dots, g_r\})\bigr),$$ \ie
		$$\textstyle[\hat U \sqcup \bigsqcup_{i \in [2,r]} \hat U \hat g_i] = [M] = [U \sqcup \bigsqcup_{i \in [2,r]} U g_i].$$
		
		Write $\hat M := \DataPhi^{(l)}\bigl((\hat U, \{\hat g_2, \dots, \hat g_r\})\bigr) = \hat U \sqcup \bigsqcup_{i \in [2,r]} \hat U \hat g_i$.
		
		Then $1 \in \hat M$ and $|[\hat M]| = |[M]| = p^{t-s+l}n$, so $\hat M \in \Omega_1^l$, so $\hat U = \Stab(\hat M)$, \cf \autoref{definition-omega}.(1) and \autoref{lemma-data-equivalences}.(2 $\Rightarrow$ 3).
		
		We have $[\hat M] = [M]$, and so $\hat M \in [M]$. Thus, we may choose $g \in G$ such that $gM = \hat M$.
		
		We have $$\hat U = \Stab(\hat M) = \Stab(gM) \stackrel{\text{R. }\ref{remark-stab-conjugation}}{=}\, {}^g\Stab(M) = {}^gU.$$
		
		So $\hat M = {}^gU \sqcup \bigsqcup_{i \in [2,r]} {}^gU \hat g_i$. Since $gM = \hat M$, we have $M = g^{-1}\hat M$, \ie
		$$U \sqcup \textstyle\bigsqcup_{i \in [2,r]} Ug_i
		= M
		= g^{-1}\hat M
		= g^{-1}({}^gU \sqcup \bigsqcup_{i \in [2,r]} {}^gU \hat g_i)
		= Ug^{-1} \sqcup \bigsqcup_{i \in [2,r]} Ug^{-1} \hat g_i.$$
		
		Again, let us write $g_1 = 1$ and $\hat g_1 := 1$ in $G$. Then
		%			$$\{Ug_1, \dots, Ug_p\} = \{Ug^{-1}\hat g_1, \dots, Ug^{-1}\hat g_p\}.$$
		$$\textstyle\bigsqcup_{j \in [1,r]} Ug_j = \textstyle\bigsqcup_{j \in [1,r]} Ug^{-1} \hat g_j.$$
		
		Then there is a unique $k \in [1,r]$ with $Ug^{-1}\hat g_1 = Ug_k$, \ie $Ug^{-1} = Ug_k$, \ie $g^{-1} \in Ug_k$, \ie $g \in g_k^{-1}U$, \ie we may choose $u \in U$ such that $g = g_k^{-1}u$. So $\hat M = gM = g_k^{-1}uM = g_k^{-1}M$ since $u \in U = \Stab(M)$. 
		
		We get $\hat U = {}^gU = {}^{g_k^{-1}u}U = {}^{g_k^{-1}}U$, and so $\{\hat g_2, \dots, \hat g_p\} \in \Tv_{r}(\hat U) = \Tv_{r}({}^{g_k^{-1}}U)$.
		
		We obtain
		$$\hat M = {}^gU \sqcup \textstyle\bigsqcup_{i \in [2,r]} {}^gU \hat g_i
		= {}^{g_k^{-1}}U \sqcup \textstyle\bigsqcup_{i \in [2,r]} {}^{g_k^{-1}}U \hat g_i$$
		and
		$$\hat M
		= g_k^{-1}M
		= g_k^{-1}(\textstyle\bigsqcup_{j \in [1,r]} Ug_j)
		= \bigsqcup_{j \in [1,r]} {}^{g_k^{-1}}Ug_k^{-1}g_j
		= {}^{g_k^{-1}}U \sqcup \bigsqcup_{j \in [1,r] \setminus \{k\}} {}^{g_k^{-1}}Ug_k^{-1}g_j.$$
		
		So $${}^{g_k^{-1}}U \sqcup \textstyle\bigsqcup_{i \in [2,r]} {}^{g_k^{-1}}U \hat g_i = {}^{g_k^{-1}}U \sqcup \bigsqcup_{j \in [1,r] \setminus \{k\}} {}^{g_k^{-1}}Ug_k^{-1}g_j.$$
		
		So $$\{{}^{g_k^{-1}}U\hat g_2, \dots, {}^{g_k^{-1}}U\hat g_r\} = \{{}^{g_k^{-1}}Ug_k^{-1}g_j : j \in [1,r] \setminus \{k\}\}.$$
		
		So $$(\hat U, \{\hat g_2, \dots, \hat g_r\}) = ({}^{g_k^{-1}}U, \{\hat g_2, \dots, \hat g_r\}) = \lambda(k) \in \lambda([1,r]).$$
		
		\emph{Injective}. Suppose given $k, k' \in [1,r]$ such that $\lambda(k) = \lambda({k'})$.
		That means $$({}^{g_k^{-1}}U, \{g_{k,2}, \dots, g_{k,r}\}) = ({}^{g_{k'}^{-1}}U, \{g_{{k'},2}, \dots, g_{{k'},r}\})$$
		with $\{g_{k,2}, \dots, g_{k,r}\} \in \Tv_{r}(^{g_k^{-1}}U)$ and $\{g_{k',2}, \dots, g_{k',r}\} \in \Tv_{r}(^{g_{k'}^{-1}}U)$ satisfying
		$$\{{}^{g_k^{-1}}Ug_{k,2}, \dots, {}^{g_k^{-1}}Ug_{k,r}\} = \{{}^{g_k^{-1}}Ug_k^{-1}g_j : j \in [1,r] \setminus \{k\}\}$$
		and
		$$\{{}^{g_{k'}^{-1}}Ug_{{k'},2}, \dots, {}^{g_{k'}^{-1}}Ug_{{k'},r}\} = \{{}^{g_{k'}^{-1}}Ug_{k'}^{-1}g_j : j \in [1,r] \setminus \{{k'}\}\}.$$
		Then ${}^{g_k^{-1}}U \sqcup \bigsqcup_{i \in [2,r]} {}^{g_k^{-1}}Ug_{k,i} = {}^{g_{k'}^{-1}}U \sqcup \bigsqcup_{i \in [2,r]} {}^{g_{k'}^{-1}}Ug_{{k'},i}$.
		
		We have \begin{align*}
			g_k^{-1}M
			&= \textstyle g_k^{-1}(\bigsqcup_{j \in [1,r]} Ug_j) \\
			&= \textstyle {}^{g_k^{-1}}U \sqcup \bigsqcup_{j \in [1,r] \setminus \{k\}} {}^{g_k^{-1}}Ug_k^{-1}g_j \\
			&= \textstyle {}^{g_k^{-1}}U \sqcup \bigsqcup_{i \in [2,r]} {}^{g_k^{-1}}Ug_{k,i} \\
			&= \textstyle {}^{g_{k'}^{-1}}U \sqcup \bigsqcup_{i \in [2,r]} {}^{g_{k'}^{-1}}Ug_{{k'},i} \\
			&= \textstyle {}^{g_{k'}^{-1}}U \sqcup \bigsqcup_{j \in [1,r] \setminus \{k'\}} {}^{g_{k'}^{-1}}Ug_{k'}^{-1}g_j \\
			&= \textstyle g_{k'}^{-1}(\bigsqcup_{j \in [1,r]} Ug_j) \\
			&= g_{k'}^{-1}M.
		\end{align*}
		
		We conclude that $g_kg_{k'}^{-1}M = M$, \ie $g_kg_{k'}^{-1} \in \Stab(M) = U$, \ie $g_k \in Ug_{k'}$, which means $Ug_k = Ug_{k'}$ and therefore $g_k = g_{k'}$, which yields $k = k'$.
		
		Hence, $\lambda$ is bijective. This shows the \emph{claim}. The assertion follows.
		
		\Ad (2). Each fibre of the map $\overline\DataPhi^{(l),l}: \Data^{(l),l} \to \overline \Omega^l$ has cardinality $p^k$ by \ref{lemma-fibres-phi-bar-1}. So \[|\Data^{(l),l}| = p^l \cdot |\overline \Omega^l|.\qedhere\]
	\end{proof}
\end{lemma}

\begin{proposition}\label{proposition-al}
	Let $l \in [0, s]$. We have $$a_l = \tfrac 1{p^l} \sum_{\substack{I \subseteq [0, l] \\ l \in I}}
	(-1)^{|I| + 1} \cdot \TvBin(\min I) \Num(s - I).$$
	\begin{proof}
		We obtain
		\begin{align*}
			a_l &\mystackrel{L. \ref{lemma-fibres-phi-bar}.(2)}{=}
			\tfrac 1{p^l} |\Data^{(l),l}| \\
			&\mystackrel{C. \ref{cor-cardinality-d-k-l}}{=}
			\tfrac 1{p^l} \sum_{\substack{I \subseteq [0, l] \\ l \in I}}
			(-1)^{|I| + 1} \cdot \TvBin(\min I) \Num(s - (I \cup \{l\})) \\
			&\mystackrel{}{=}
			\tfrac 1{p^l} \sum_{\substack{I \subseteq [0, l] \\ l \in I}}
			(-1)^{|I| + 1} \cdot \TvBin(\min I) \Num(s - I). \qedhere
		\end{align*}
	\end{proof}
\end{proposition}

%% file: text/BA-sub-formula-a_l.tex
\begin{theorem}\label{thm-sub-formel-pl}
	Recall that $G$ is a finite group of order $|G| = p^tn$.
	
	Let $s \in [0, t]$. Let $l \in [0, s]$.
	
	Recall that, given $I = \{c_1, \dots, c_k\} \subseteq [0, l]$ where $c_1 < \ldots < c_k$, we denote by $\Num(s-I)$ the number of chains of $p$-subgroups $U_1 \leq \ldots \leq U_k \leq G$, where $U_1 \in \Sub(s-c_1), \dots, U_k \in \Sub(s-c_k)$, \cf \autoref{definition-sub-num}\upshape.(2).
	
	\itshape We have the \emph{sieve formula} $$\sum_{\substack{I \subseteq [0, l] \\ I \neq \emptyset}} (-1)^{|I| + 1} \Num(s - I) \equiv_{p^{l+1}} 1.$$
	\begin{proof}
		\emph{Claim}. We have $$\sum_{\substack{I \subseteq [0, l] \\ I \neq \emptyset}} (-1)^{|I|+1}(\Num(s-I)-1) \equiv_{p^{l+1}} 0.$$
		
		Write $f(J) :=(-1)^{|I| + 1} (\Num(J) - 1) \in \Z$ for $J \subseteq [s-l, s]$.
		
		The claim amounts to showing $$\sum_{\substack{I \subseteq [0, l] \\ I \neq \emptyset}} f(s-I) \stackrel!{\equiv}_{p^{l+1}} 0. \leqno{(\Diamond_{l,s})}$$
		
		Let $X := \{(l, s) \in \Z_{\geq 0} \times \Z_{\geq 0}: l \leq s \leq t\}$. For $(l, s), (l', s') \in X$, we define: $$(l, s) \preccurlyeq (l', s') \ :\Leftrightarrow\ l \leq l' \text{ and } s \leq s'.$$
		Then $(X, \preccurlyeq)$ is a partially ordered set. By \autoref{cor-finite-posets-noetherian}, $(X, \preccurlyeq)$ is noetherian.
		
		Therefore, we can use Noetherian induction over $(X, \preccurlyeq)$ to prove $(\Diamond_{l,s})$ for $(l, s) \in X$, \cf \autoref{lemma-noetherian-induction}.
		
		Suppose given $(l, s) \in X$. We may suppose that $$\sum_{\substack{I \subseteq [0, l'] \\ I \neq \emptyset}} f(s'-I) \equiv_{p^{l'+1}} 0 \leqno{(\Diamond_{l',s'})}$$holds for $(l', s') \in X$ with $(l', s') \prec (l, s)$. We have to show $(\Diamond_{l,s})$.
		
%		\textit{Base case}. For $l = 0$, we get $$\sum_{I \subseteq [0, 0]} f(s-I) = \sum_{I \subseteq [0, 0]} (-1)^{|I| + 1} (\Num(s - I) - 1) = \Num(s) - 1 \equiv_p 0$$ by \autoref{thm-sylow-wielandt}. Note that $\Num(s - \emptyset) = \Num(\emptyset) = |\{()\}| = 1$.
		
		\textit{Induction step}. By \autoref{remark-q-mod-p^k}, we have $$q \equiv_{p^{l+1}} \sum_{k \in [0, l]} a_kp^k$$ where
		$$a_k = \tfrac 1{p^k} \sum_{\substack{I \subseteq [0, k] \\ k \in I}}
		(-1)^{|I| + 1} \cdot \TvBin(\min I) \Num(s - I)$$ for $k \in [0, l]$, thanks to \autoref{proposition-al}. Then
		\[\def\arraystretch{2}
			\begin{array}{cll}
				q &\equiv_{p^{l+1}} &\displaystyle\sum_{k \in [0, l]} a_kp^k \\
				&= &\displaystyle\sum_{k \in [0, l]} \displaystyle\sum_{\substack{I \subseteq [0, k] \\ k \in I}}
				(-1)^{|I| + 1} \cdot \TvBin(\min I) \Num(s - I) \\
				&= &\displaystyle\sum_{\substack{I \subseteq [0, l] \\ I \neq \emptyset}}
				(-1)^{|I| + 1} \cdot \TvBin(\min I) \Num(s - I) \\
				&= &\displaystyle\sum_{k \in [0, l]} \displaystyle\sum_{\substack{I \subseteq [k, l] \\ k \in I}}
				(-1)^{|I| + 1} \cdot \TvBin(\min I) \Num(s - I) \\
				&= &\displaystyle\sum_{k \in [0, l]} \TvBin(k) \displaystyle\sum_{\substack{I \subseteq [k, l] \\ k \in I}}
				(-1)^{|I| + 1} \Num(s - I).
			\end{array}
		\]
		Considering the cyclic group $C := \Cycl_{p^tn}$ with $p^tn$ elements, we have $\Num_{C,p}(s-I) = 1$ for $I \subseteq [0, l]$ and therefore $$q \equiv_{p^{l+1}} \sum_{k \in [0, l]} \TvBin(k) \sum_{\substack{I \subseteq [k, l] \\ k \in I}}
		(-1)^{|I| + 1}.$$
		So
		\[\def\arraystretch{2}
			\begin{array}{clll}
				0 &\equiv_{p^{l+1}} &&\displaystyle\Biggl(\sum_{k \in [0, l]} \TvBin(k) \displaystyle\sum_{\substack{I \subseteq [k, l] \\ k \in I}}
				(-1)^{|I| + 1} \Num(s - I)\Biggr) - q \\
				&\equiv_{p^{l+1}} &&\displaystyle\Biggl(\sum_{k \in [0, l]} \TvBin(k) \displaystyle\sum_{\substack{I \subseteq [k, l] \\ k \in I}}
				(-1)^{|I| + 1} \Num(s - I)\Biggr) - \Biggl(\sum_{k \in [0, l]} \TvBin(k) \displaystyle\sum_{\substack{I \subseteq [k, l] \\ k \in I}}
				(-1)^{|I| + 1}\Biggr) \\
				&=&&\displaystyle\sum_{k \in [0, l]} \TvBin(k) \displaystyle\sum_{\substack{I \subseteq [k, l] \\ k \in I}}
				(-1)^{|I| + 1} (\Num(s - I) - 1) \\
				&=&&\displaystyle\sum_{k \in [0, l]} \TvBin(k) \displaystyle\sum_{\substack{I \subseteq [k, l] \\ k \in I}}
				f(s-I).
			\end{array}
		\]
		\textsf{Subclaim}. We have $$\sum_{k \in [0, l]} \TvBin(k) \sum_{\substack{I \subseteq [k, l] \\ k \in I}} f(s-I) \equiv_{p^{l+1}} \Biggl(\sum_{k \in [0, l-j-1]} \TvBin(k) \sum_{\substack{I \subseteq [k, l] \\ k \in I}}
		f(s-I)\Biggr) + \TvBin(l-j)\Biggl(\sum_{\substack{I \subseteq [l-j, l] \\ I \neq \emptyset}}
		f(s-I)\Biggr).$$
		for $j \in [0, l]$.
		
		To prove the subclaim, we proceed by induction on $j$. For the base case $j = 0$, we obtain
		$$\sum_{k \in [0, l]} \TvBin(k) \sum_{\substack{I \subseteq [k, l] \\ k \in I}} f(s-I) = \Biggl(\sum_{k \in [0, l-1]} \TvBin(k) \sum_{\substack{I \subseteq [k, l] \\ k \in I}}
		f(s-I)\Biggr) + \TvBin(l)\Biggl(\sum_{\substack{I \subseteq [l, l] \\ I \neq \emptyset}}
		f(s-I)\Biggr).$$
		
		Let $j \in [0, l-1]$. Suppose that the statement holds for $j$. We show the statement for $j+1$.
		
		We see that
		\[\def\arraystretch{2}
			\begin{array}{ll}
				&\displaystyle\sum_{k \in [0, l]} \TvBin(k) \sum_{\substack{I \subseteq [k, l] \\ k \in I}} f(s-I) \\
				\stackrel{\textsf{IH}}{\equiv}_{p^{l+1}}
				&\Biggl(\displaystyle\sum_{k \in [0, l-j-1]} \TvBin(k) \sum_{\substack{I \subseteq [k, l] \\ k \in I}}
				f(s-I)\Biggr) + \TvBin(l-j)\Biggl(\sum_{\substack{I \subseteq [l-j, l] \\ I \neq \emptyset}} f(s-I)\Biggr)\\
				=
				&\Biggl(\displaystyle\sum_{k \in [0, l-j-1]} \TvBin(k) \sum_{\substack{I \subseteq [k, l] \\ k \in I}}
				f(s-I)\Biggr) + \TvBin(l-j)\Biggl(\sum_{\substack{I \subseteq [0, j] \\ I \neq \emptyset}}
				f((s-(l-j))-I)\Biggr).
			\end{array}
		\]
		We have $(j, s-(l-j)) \prec(l, s)$. Since $$\sum_{\substack{I \subseteq [0, j] \\ I \neq \emptyset}}
		f((s-(l-j))-I) \equiv_{p^{j+1}} 0 \leqno{(\Diamond_{j, s-(l-j)})}$$ holds thanks to our outer induction hypothesis of the proof of the claim, and since
		$$\TvBin(l-j) \equiv_{p^{l-j}} \TvBin(l-j-1)$$ thanks to \autoref{lemma-b-k-mod-p-k}, we conclude that
		$$\TvBin(l-j)\Biggl(\sum_{\substack{I \subseteq [0, j] \\ I \neq \emptyset}}
		f((s-(l-j))-I)\Biggr)
		\equiv_{p^{l+1}}
		\TvBin(l-j-1)\Biggl(\sum_{\substack{I \subseteq [0, j] \\ I \neq \emptyset}}
		f((s-(l-j))-I)\Biggr).$$
		\allowdisplaybreaks
		So we may continue to get
		\[\def\arraystretch{2}
			\begin{array}{lcl}
				&&\Biggl(\displaystyle\sum_{k \in [0, l-j-1]} \TvBin(k) \sum_{\substack{I \subseteq [k, l] \\ k \in I}}
				f(s-I)\Biggr) + \TvBin(l-j)\Biggl(\sum_{\substack{I \subseteq [0, j] \\ I \neq \emptyset}}
				f((s-(l-j))-I)\Biggr) \\
				\equiv_{p^{l+1}}&&\Biggl(\displaystyle\sum_{k \in [0, l-j-1]} \TvBin(k) \sum_{\substack{I \subseteq [k, l] \\ k \in I}}
				f(s-I)\Biggr) + \TvBin(l-j-1)\Biggl(\sum_{\substack{I \subseteq [0, j] \\ I \neq \emptyset}}
				f((s-(l-j))-I)\Biggr) \\
				=&&\displaystyle\Biggl(\sum_{k \in [0, l-j-2]} \TvBin(k) \sum_{\substack{I \subseteq [k, l] \\ k \in I}} f(s-I)\Biggr)\\
				&+&\displaystyle\TvBin(l-j-1)\Biggl(\sum_{\substack{I \subseteq [l-j-1, l] \\ l-j-1 \in I}} f(s-I) + \sum_{\substack{I \subseteq [0, j] \\ I \neq \emptyset}}
				f((s-(l-j))-I)\Biggr) \\
				=&&\displaystyle\Biggl(\sum_{k \in [0, l-j-2]} \TvBin(k) \sum_{\substack{I \subseteq [k, l] \\ k \in I}} f(s-I)\Biggr)\\
				&+&\displaystyle\TvBin(l-j-1)\Biggl(\sum_{\substack{I \subseteq [l-j-1, l] \\ l-j-1 \in I}} f(s-I) + \sum_{\substack{I \subseteq [l-j, l] \\ I \neq \emptyset}}
				f(s-I)\Biggr) \\
				=&&\Biggl(\displaystyle\sum_{k \in [0, l-j-2]} \TvBin(k) \sum_{\substack{I \subseteq [k, l] \\ k \in I}} f(s-I)\Biggr) + \TvBin(l-j-1)\Biggl(\sum_{\substack{I \subseteq [l-j-1, l] \\ I \neq \emptyset}} f(s-I)\Biggr).
			\end{array}
		\]
		
		This shows the \textsf{subclaim}. Using the subclaim for $j = l$ yields
		\[\def\arraystretch{2}
			\begin{array}{clll}
				0 &\equiv_{p^{l+1}} &&\displaystyle\sum_{k \in [0, l]} \TvBin(k) \displaystyle\sum_{\substack{I \subseteq [k, l] \\ k \in I}}
				f(s-I) \\
				&\equiv_{p^{l+1}} &&\displaystyle\Biggl(\sum_{k \in [0, -1]} \TvBin(k) \sum_{\substack{I \subseteq [k, l] \\ k \in I}}
				f(s-I)\Biggr) + \TvBin(0)\Biggl(\sum_{\substack{I \subseteq [0, l] \\ I \neq \emptyset}}
				f(s-I)\Biggr) \\
				&= &&\displaystyle\sum_{\substack{I \subseteq [0, l] \\ I \neq \emptyset}}
				f(s-I) \\
				&= &&\displaystyle\sum_{\substack{I \subseteq [0, l] \\ I \neq \emptyset}}
				(-1)^{|I| + 1} (\Num(s - I) - 1).
			\end{array}
		\]
		This shows the \emph{claim}.
		
		To simplify the term further, note that
		$$\sum_{\substack{I \subseteq [0, l] \\ I \neq \emptyset}} (-1)^{|I|+1} = 1 - \sum_{I \subseteq [0, l]}
		(-1)^{|I|} = 1 - \sum_{m \in [0, l]} (-1)^m\binom{l}{m} = 1.$$
		So $$0 \equiv_{p^{l+1}} \sum_{\substack{I \subseteq [0, l] \\ I \neq \emptyset}}
		(-1)^{|I| + 1} (\Num(s - I) - 1) = \Biggl(\sum_{\substack{I \subseteq [0, l] \\ I \neq \emptyset}}
		(-1)^{|I| + 1} \Num(s - I)\Biggr) - 1,$$
		
		\ie \[\sum_{\substack{I \subseteq [0, l] \\ I \neq \emptyset}}
		(-1)^{|I| + 1} \Num(s - I) \equiv_{p^{l+1}} 1.\qedhere\]
	\end{proof}
\end{theorem}

\begin{remark without proof}
	\autoref{thm-sub-formel-pl} has formal similarities with the sieve formula from set theory:
	
	Let $X$ be a finite set and $Y_0, \dots, Y_l \subseteq X$ be subsets such that $X = \bigcup_{k \in [0,l]} Y_k$. Let $$Y_I := \textstyle\bigcap_{i \in I} Y_i$$ for $I \subseteq [0,l]$. Then $$\sum_{\substack{I \subseteq [0,l] \\ I \neq \emptyset}} (-1)^{|I|+1} |Y_I| = |X|.$$
\end{remark without proof}

\begin{remark without proof}
	Let $s \in [0, t]$. Consider the case $l = s$ in \autoref{thm-sub-formel-pl}. We have
	\begin{align*}
		&\sum_{I \subseteq [0, s]} (-1)^{|I| + 1} \Num(s - I) \\
		\mystackrel{}{=} &\sum_{\substack{I \subseteq [0, s] \\ s \in I}} (-1)^{|I| + 1} \Num(s - I) + \sum_{\substack{I \subseteq [0, s] \\ s \notin I}} (-1)^{|I| + 1} \Num(s - I) \\
		\mystackrel{}{=} &\sum_{I \subseteq [0, s-1]} (-1)^{|I| + 1} \Num(s - (I \sqcup \{s\})) + \sum_{I \subseteq [0, s-1]} (-1)^{|I|} \Num(s - I) \\
		\mystackrel{R. \ref{remark-trivial-subgroup-in-chain}}{=} &\sum_{I \subseteq [0, s-1]} (-1)^{|I| + 1} \Num(s - I) + \sum_{I \subseteq [0, s-1]} (-1)^{|I|} \Num(s - I) \\
		\mystackrel{}{=} &0,
	\end{align*} \ie $$\sum_{\substack{I \subseteq [0, s] \\ I \neq \emptyset}} (-1)^{|I| + 1} \Num(s - I) = 1.$$
\end{remark without proof}

%% file: text/BA-coherence-Sylow.tex
\begin{remark without proof}\label{remark-abkürzung-sylow}
	One could ask whether the formula from \autoref{thm-sub-formel-pl} follows directly by the Theorem of Sylow-Frobenius.
	
	We have $$\Num(s-l) - 1 \equiv_p 0$$ for $l \in [0, s]$, \cf \autoref{thm-sylow-wielandt},
	and therefore $$\prod_{k \in [0, l]} (\Num(s-k) - 1) = \sum_{I \subseteq [0, l]} (-1)^{|I|}\prod_{i \in I} \Num(s-i) \equiv_{p^{l+1}} 0$$ for $l \in [0,s]$.
	
	We consider some cases.
	\begin{enumerate}[listparindent=0pt]
		\item Let $l = 1$. Then
		\begin{center}
			$\begin{aligned}
				&&&\textstyle\sum_{I \subseteq [0, 1]} (-1)^{|I|}\prod_{i \in I} \Num(s-i) \\
				&= &&(-1)^0 \cdot 1 + (-1)^1(\Num(s) + \Num(s{-1})) + (-1)^2 \Num(s{-1})\Num(s) \\
				&\equiv_{p^2} &&0,
			\end{aligned}$
		\end{center} \ie $$\Num(s) + \Num(s{-1}) - \Num(s{-1})\Num(s) \equiv_{p^2} 1.$$
		By \autoref{thm-sub-formel-pl}, we have $$\Num(s) + \Num(s{-1}) - \Num(s{-1}, s) \equiv_{p^2} 1.$$
		Hence $$\Num(s{-1})\Num(s) \equiv_{p^2} \Num(s{-1}, s).$$
		\item Let $l = 2$. Then
		\begin{center}
			$\begin{aligned}
				&&&\textstyle\textover{}{\sum_{I \subseteq [0, 2]} (-1)^{|I|}\prod_{i \in I} \Num(s-i)} \\
				&= &&(-1)^0 \cdot 1 \hspace{-10pt} &&+ (-1)^1(\textover{}{\Num(s) + \Num(s{-1}) + \Num(s{-2}))} \\
				&&& &&+ (-1)^2 (\Num(s{-1})\Num(s) + \Num(s{-2})\Num(s) + \Num(s{-2})\Num(s{-1})) \\
				&&& &&+ (-1)^3 \textover{}{\Num(s{-2})\Num(s{-1})\Num(s)} \\
				&\equiv_{p^3} &&0,
			\end{aligned}$
		\end{center} \ie
		\begin{align*}
			\Num(s) &+ \Num(s{-1}) + \Num(s{-2}) - \Num(s{-1})\Num(s) \\
			&- \Num(s{-2})\Num(s) - \Num(s{-2})\Num(s{-1}) + \Num(s{-2})\Num(s{-1})\Num(s) \equiv_{p^3} 1.
		\end{align*}
		By \autoref{thm-sub-formel-pl}, we have 
		\begin{align*}
			\Num(s) &+ \Num(s{-1}) + \Num(s{-2}) - \Num(s{-1},s) \\
			&- \Num(s{-2},s) - \Num(s{-2}, s{-1}) + \Num(s{-2}, s{-1}, s) \equiv_{p^3} 1.
		\end{align*}
		Hence
		\begin{align*}
			&\phantom{\equiv_{p^3}}\ \Num(s{-1})\Num(s) + \Num(s{-2})\Num(s) + \Num(s{-2})\Num(s{-1}) - \Num(s{-2})\Num(s{-1})\Num(s) \\
			&{\equiv_{p^3}}\ \Num(s{-1},s) + \Num(s{-2},s) + \Num(s{-2}, s{-1}) - \Num(s{-2}, s{-1}, s).
		\end{align*}
		This congruence does not necessarily imply that
		$$\Num(s{-2})\Num(s{-1})\Num(s) \equiv_{p^3} \Num(s{-2}, s{-1}, s).$$
		For instance, if we take $G = \Sym_5$, $p = 2$ and $s = 3$, then Magma gives $$\Num_{\Sym_5, 2}(1) = 25, \quad \Num_{\Sym_5, 2}(2) = 35, \quad \Num_{\Sym_5, 2}(3) = 15,$$ and
		$$\Num_{\Sym_5, 2}(1, 2, 3) = 105,$$ but $25 \cdot 35 \cdot 15 = 13125 \equiv_8 5 \not \equiv_8 1 \equiv_8 105$.
		
		\item The example in (2) shows that we do not have the congruence $$\prod_{i \in I} \Num(i) \equiv_{p^{|I|}} \Num(I)$$ for $\emptyset \neq I \subseteq [0, t]$ in general.
		
		So it seems to be impossible to conclude \autoref{thm-sub-formel-pl} from \autoref{thm-sylow-wielandt} solely.
	\end{enumerate}
	
\end{remark without proof}

%% file: text/BA-sub-formula-a_1.tex
\ifarxiv\else
	\begin{remark}\label{remark-a_1-siebformel-direkt}
		We have $$\Num(s) + \Num(s-1) - \Num(s-1,s) \equiv_{p^2} 1.$$
		\begin{proof}
			By \autoref{remark-q-mod-p^k}, we have $$q \equiv_{p^2} a_0 + a_1p,$$ where
			$$a_0 = \Num(s),$$ and $$a_1 = \tfrac 1p (\TvBin(1)\Num(s-1) - \Num(s-1, s)),$$ thanks to Remarks \ref{remark-a_0-direkt} and \ref{remark-a_1-direkt}.
	
			Then
			\begin{align*}
				q\quad
				\stackrel{}{\equiv_{p^2}} \quad&a_0 + a_1p \\
				\stackrel{}{\textover{\equiv_{p^2}}{=}} \quad&\Num(s) + \TvBin(1)\Num(s-1) - \Num(s-1, s).
			\end{align*}
			Considering $C := \C_{p^tn}$, we get $$q \equiv_{p^2} \Num_{C,p}(s) + \TvBin(1)\Num_{C,p}(s-1) - \Num_{C,p}(s-1, s) = 1 + \TvBin(1) - 1.$$
			Since $\TvBin(1) \equiv_{p} \TvBin(0) = 1$ by \autoref{lemma-b-k-mod-p-k} and $\Num(s-1) \equiv_p 1$ by \autoref{thm-sylow-wielandt}, we conclude that
			\begin{align*}
				0\quad &\equiv_{p^2}\quad (\Num(s) - 1) + \TvBin(1)(\Num(s-1) - 1) - (\Num(s-1, s) - 1) \\
				&\equiv_{p^2}\quad (\Num(s) - 1) + 1 \cdot (\Num(s-1) - 1) - (\Num(s-1, s) - 1) \\
				&\equiv_{p^2}\quad \Num(s) - 1 + \Num(s-1) - \Num(s-1, s).\qedhere
			\end{align*}
			Note that e.g. for $p = 2$ and $s \in [1, t-1]$, we get $$\TvBin(1) = \binom{2^{t-s+1}n - 1}{2-1} = 2^{t-s+1}n - 1 \equiv_{2^2} 3.$$
		\end{proof} % geprüft bis Gruppen der Ordnung 255.
	\end{remark}
\fi

\begin{example}\label{ex-p^2-formula}
	The following examples have been obtained using the computer algebra system Magma \bfcit{Mag}.
	
	\begin{enumerate}[listparindent=0pt]
		\item Let $G = \Sym_3$ and $p = 2$. Magma gives the following.
		$$\Num(1) = 3, \qquad \Num(0) = 1, \qquad \Num(0, 1) = 3.$$
		Then we get the following result for $s = 1$.
		$$\Num(1) + \Num(0) - \Num(0, 1) = 3 + 1 - 3 = 1 \equiv_4 1,$$ \Cf also \autoref{remark-trivial-subgroup-in-chain}.
		
		\item Let $G = \Sym_4$ and $p = 2$. Magma gives the following.
		\begin{align*}
			\Num(3) &= 3, &
			\Num(2) &= \textover{15}{7},&
			\Num(1) &= 9,&&
			\Num(0) = 1, \\
			\Num(2, 3) &= 9,&
			\Num(1, 2) &= 15,&
			\Num(0, 1) &= 9.
		\end{align*}
		
		\Cf also \autoref{ex-sym-4-sub-num}.
		
		Then we get the following results for $s = 3$, $s = 2$ and $s = 1$, respectively.
		\[
		\begin{array}{lllllllllllrlllll}
			&\Num(3) &+ &\Num(2) &- &\Num(2, 3) &= &3 &+ &7 &- &9  &= 1 &\equiv_4 1, \\
			&\Num(2) &+ &\Num(1) &- &\Num(1, 2) &= &7 &+ &9 &- &15 &= 1 &\equiv_4 1, \\
			&\Num(1) &+ &\Num(0) &- &\Num(0, 1) &= &9 &+ &1 &- &9  &= 1 &\equiv_4 1.
		\end{array}
		\]
		
		\item Let $G = \Sym_6$ and $p = 3$. Then $|G| = 720 = 3^2 \cdot 80$. Therefore, we have $t = 2$ and $n = 80$.
		
		We choose $s = 2 \in [1,2]$. Magma gives the following.
		$$\Num_{\Sym_6, 3}(2) = 10, \qquad \Num_{\Sym_6, 3}(1)= 40, \qquad \Num_{\Sym_6, 3}(1, 2) = 40.$$
		Then
		$$\Num_{\Sym_6, 3}(2) + \Num_{\Sym_6, 3}(1) - \Num_{\Sym_6, 3}(1, 2) = 10 + 40 - 40 = 10 \equiv_9 1.$$
		Note that $10 \not \equiv_{27} 1$.
		
		\item We have $$\Num_{\Alt_5, 2}(2) + \Num_{\Alt_5, 2}(1) - \Num_{\Alt_5, 2}(1,2)
		= 5 + 15 - 15 = 5 \equiv_4 1.$$
		Note that $5 \not \equiv_8 1$.
		
		\item We have $$\Num_{\Alt_6, 2}(3) + \Num_{\Alt_5, 6}(2) - \Num_{\Alt_6, 2}(2,3)
		= 45 + 75 - 135 = -15 \equiv_4 1.$$
		
		\item We have $$\Num_{\GL_3(\F_2), 2}(2) + \Num_{\GL_3(\F_2), 2}(1) - \Num_{\GL_3(\F_2), 2}(1, 2)
		= 35 + 21 - 63 = -7 \equiv_4 1.$$
		
		\item We have $$\Num_{\GL_3(\F_5), 2}(5) + \Num_{\GL_3(\F_5), 2}(4) - \Num_{\GL_3(\F_5), 2}(4, 5) =  15 +  35 - 45 =  5 \equiv_4 1.$$
		
		\item We have $$\Num_{\SL_3(\F_3), 2}(2) + \Num_{\SL_3(\F_3), 2}(1) - \Num_{\SL_3(\F_3), 2}(1, 2) = 585 + 117 - 1053 = -351 \equiv_4 1.$$
		
		\item \newcommand{\Math}{\operatorname{M}}
		We consider the Mathieu group $$\Math_{11} := \langle (1,2,3,4,5,6,7,8,9,10,11), (3,7,11,8)(4,10,5,6)\rangle \leq \Sym_{11}.$$ Note that $|\Math_{11}| = 7920 = 2^4 \cdot 495$. We have
		$$\Num_{\Math_{11}, 2}( 4) + \Num_{\Math_{11}, 2}( 3) - \Num_{\Math_{11}, 2}( 3,  4) = 495 + 1155 -1485 = 165 \equiv_4 1.$$
	\end{enumerate}
\end{example}

%% file: text/BA-sub-formula-a_2.tex

\begin{example}
	Let $$G = \textrm F_9 = \spn{a,b,c \mid a^3,\, b^3,\, c^3,\, [a,b],\, a^c=b, b^c = ab^2} \backsimeq (\Cycl_3 \times \Cycl_3) \rtimes \Cycl_8$$ be the Frobenius group of order $72$.
	
	Let $p = 2$. Then Magma produces the following result for $s = 3$.
	 \[
	 	\def\arraystretch{1.25}
	 	\begin{array}{cl}
	 		
	 		&\Num(3) + \Num(2) + \Num(1) - \Num(2,3) - \Num(1,3) - \Num(1,2) + \Num(1,2,3) \\
	 		\stackrel{}{\textover{\equiv_8}{=}}\ &9 + 9 + 9 - 9 - 9 - 9 + 9 \\
	 		\stackrel{}{\textover{\equiv_8}{=}}\ &9 \\
	 		\stackrel{}{\equiv_8} &1.
	 	\end{array}
	\] Note that $9 \not \equiv_{16} 1$.
\end{example}

%% file: text/BA-questions.tex
\begin{question}
	Let $s \in [0, t]$. Let $J \subseteq [0, s]$ such that $J \neq \emptyset$.
	
	Does the congruence $$\sum_{\substack{I \subseteq J \\ I \neq \emptyset}} (-1)^{|I| + 1} \Num(s - I) \equiv_{p^{|J|}} 1$$ hold?
\end{question}

\begin{question}\label{qst-p-group-congruence}
	Let $s \in [0, t]$. Suppose $G$ to be a $p$-group, \ie $|G| = p^t$. Let $l \in [0, s]$.
	
	Does the congruence $$\sum_{\substack{I \subseteq [0, l] \\ I \neq \emptyset}} (-1)^{|I| + 1} \Num(s - I) \equiv_{p^{\frac{(l+1)(l+2)}{2}}} 1$$ hold?
\end{question}

\begin{remark without proof}
	Investigating \autoref{qst-p-group-congruence} via Magma \bfcit{Mag}, where we used the Small Groups Library \bfcit{Bes}, gives the following results.
	\begin{enumerate}[listparindent=0pt]
		\item Let $s \in [1, t]$. Let $l = 1$. Then
		$$\Num(s) + \Num(s-1) - \Num(s-1,s) \equiv_{p^3} 1$$
		holds for $p$-groups of order $$|G| \in \{2^1,\dots,2^7, 3^1,\dots, 3^5, 5^1, \dots, 5^5, 7^1, \dots, 7^3\}.$$
		
		\item Let $s \in [2, t]$. Let $l = 2$. Then
		\begin{align*}
			1\ \equiv_{p^6}\ \phantom{+\,} &\Num(s) + \Num(s{-1}) + \Num(s{-2}) \\
			-\, &\Num(s{-1},s) - \Num(s{-2},s) - \Num(s{-2}, s{-1}) \\
			+\, &\Num(s{-2}, s{-1}, s)
		\end{align*}
		holds for $p$-groups of order $$|G| \in \{2^2,\dots,2^7, 3^2,\dots, 3^5, 5^2, \dots, 5^4, 7^2, 7^3\}.$$
		
		\item Let $s \in [3, t]$. Let $l = 3$. Then
		\begin{align*}
			1\ \equiv_{p^{10}}\ \phantom{+\,} &\Num(s) + \Num(s{-1}) + \Num(s{-2}) + \Num(s{-3}) \\
			-\, &\Num(s{-1},s) - \Num(s{-2},s) - \Num(s{-2},s{-1}) \\
			-\, &\Num(s{-3}, s) - \Num(s{-3}, s{-1}) - \Num(s{-3},s{-2}) \\
			+\, &\Num(s{-2}, s{-1}, s) + \Num(s{-3}, s{-1}, s) + \Num(s{-3}, s{-2}, s) + \Num(s{-3}, s{-2}, s{-1}) \\
			-\, &\Num(s{-3}, s{-2}, s{-1}, s)
		\end{align*}
		holds for $p$-groups of order $$|G| \in \{2^3,\dots,2^6, 3^3,\dots, 3^5, 5^3, 7^3\}.$$
		
		\item Let $s \in [4, t]$. Let $l = 4$. Then
		\begin{align*}
			1\ \equiv_{p^{15}}\ \phantom{+\,} &\Num(s) + \Num(s{-1}) + \Num(s{-2}) + \Num(s{-3}) + \Num(s{-4}) \\
			-\, &\Num(s{-1},s) - \Num(s{-2},s) - \Num(s{-2},s{-1}) \\
			-\, &\Num(s{-3}, s) - \Num(s{-3}, s{-1}) - \Num(s{-3},s{-2}) \\
			-\, &\Num(s{-4}, s) - \Num(s{-4}, s{-1}) - \Num(s{-4},s{-2}) - \Num(s{-4},s{-3}) \\
			+\, &\Num(s{-2}, s{-1}, s) + \Num(s{-3}, s{-1}, s) + \Num(s{-3}, s{-2}, s) \\
			+\, &\Num(s{-3}, s{-2}, s{-1}) + \Num(s{-4}, s{-1}, s) + \Num(s{-4}, s{-2}, s) \\
			+\, &\Num(s{-4}, s{-2}, s{-1}) + \Num(s{-4}, s{-3}, s) + \Num(s{-4}, s{-3}, s{-1}) \\
			+\, &\Num(s{-4}, s{-3}, s{-2}) \\
			-\, &\Num(s{-3}, s{-2}, s{-1}, s) - \Num(s{-4}, s{-2}, s{-1}, s) - \Num(s{-4}, s{-3}, s{-1}, s) \\
			-\, &\Num(s{-4}, s{-3}, s{-2}, s) - \Num(s{-4}, s{-3}, s{-2}, s{-1}) \\
			+\, &\Num(s{-4}, s{-3}, s{-2}, s{-1}, s)
		\end{align*}
		holds for $p$-groups of order $$|G| \in \{2^4, \dots, 2^6, 3^4, 3^5\}.$$
	\end{enumerate}
\end{remark without proof}

\begin{remark without proof}
	Examples \ref{ex-p^2-formula}.(3, 4) show that the assumption that $G$ is a $p$-group in \autoref{qst-p-group-congruence} cannot be omitted.
\end{remark without proof}

%% file: text/BA-formalia.tex
% Offizieller Teil als Anhang
\selectlanguage{german}

\addchap*{Zusammenfassung}
\pagenumbering{gobble}
\input{text/BA-summary-German.tex}

\addchap*{Versicherung}
Hiermit versichere ich, Elias Schwesig,
\begin{enumerate}[listparindent=0pt]
	\item dass ich meine Arbeit selbstständig verfasst habe,
	\item dass ich keine anderen als die angegebenen Quellen benutzt und alle wörtlich oder
	sinngemäß aus anderen Werken übernommenen Aussagen als solche
	gekennzeichnet habe,
	\item dass die eingereichte Arbeit weder vollständig noch in wesentlichen Teilen
	Gegenstand eines anderen Prüfungsverfahrens gewesen ist,
	\item dass das elektronische Exemplar mit den anderen Exemplaren übereinstimmt. 
\end{enumerate}

\begin{flushright}
	Stuttgart, im Juli 2024
	
	\vspace{1cm}
	Elias Schwesig
\end{flushright}

%% file: text/BA-summary-German.tex
% TODO Hyphenation auf Deutsch umstellen?
Sei $G$ eine endliche Gruppe. Sei $p$ eine Primzahl.

Wir schreiben $|G| = p^tn$, wobei $t = \val_p(|G|) \in \Z_{\geq 0}$, $n \in \Z_{\geq 1}$ und $n \not \equiv_p 0$.

Sei $s \in [0, t]$.

Für $k \geq 0$ und $I = \{c_1, \dots, c_k\} \subseteq [0, t]$, wobei $c_1 > \ldots > c_k$, schreiben wir $$\Num(s-I) := |\{(U_1, \dots, U_k) : U_1 \leq U_2 \leq \ldots \leq U_k \leq G,\ |U_i| = p^{s-c_i} \text{ für } i \in [1,k]\}|.$$ Es bezeichnet also $\Num(s-I)$ die Anzahl der Ketten von $p$-Untergruppen $$U_1 \leq U_2 \leq \ldots \leq U_k \leq G,$$ wobei $U_i$ Ordnung $p^{s-c_i}$ hat für $i \in [1,k]$.

Wir erhalten folgende Siebformel für $l \in [0, s]$.
\begin{align*}
	\sum_{\substack{I \subseteq [0, l] \\ I \neq \emptyset}} (-1)^{|I| + 1} \Num(s-I) \equiv_{p^{l+1}} 1.
\end{align*}
Der Name \glqq Siebformel\grqq\ entlehnt sich der Tatsache, dass die Kongruenz formale Ähnlichkeiten mit der Siebformel aus der Mengenlehre hat, die auch als Prinzip von Inklusion und Exklusion bekannt ist.

Für $l = 0$ wurde die Siebformel im Fall $s = t$ bereits von Sylow \bfcite{Syl}{Th.\,II} gezeigt und als ein Teil der Sylowsätze bekannt. Frobenius \bfcite{Fro}{§4, I.} hat sie für den allgemeinen Fall $s \in [0, t]$ bewiesen.

Im Jahr 1959 hat Wielandt \bfcit{Wie} einen alternativen Beweis im Fall $l = 0$ angeführt, in dem er eine Gruppenoperation von $G$ auf der Menge aller Teilmengen von $G$, die $p^s$ Elemente haben, definiert, und dann durch geschicktes Abzählen von Bahnen gewisser Länge die Kongruenz $$\binom{p^tn}{p^s} \equiv_p \Num(s- \{0\})$$ erhält, die nur von $|G|$ und nicht von $G$ selbst abhängt.

Graham Higman verkürzte Wielandts Beweis, indem er diese Kongruenz auf die zyklische Gruppe $\Cycl_{p^tn}$ angewandt hat und so ohne weitere zahlentheoretische Überlegungen die Kongruenz $$\binom{p^tn}{p^s} \equiv_p 1$$ erhielt. Es folgt in diesem Fall die behauptete Aussage $$\Num(s - \{0\}) \equiv_p 1.$$

Beim Beweis der Siebformel sind wir Wielandts Idee gefolgt, wir haben uns jedoch nicht nur auf das Abzählen von Bahnen gewisser Länge beschränkt, sondern haben Bahnen beliebiger Länge abgezählt. Hat man diese Anzahlen in der Hand, ermöglicht das einem, eine Kongruenz modulo $p^{l+1}$ aufzustellen, und in einem weiteren Schritt zur behaupteten Siebformel umzuformen, wobei sich wie bei Higman ein Vergleich mit dem Fall der zyklischen Gruppe als hilfreich erwies.

Wir haben zudem Beispiele angeführt, die im Fall $l = 1$ und $l = 2$ zeigen, dass die Kongruenz nicht verbessert werden kann.

%% file: text/BA-a_0.tex
Let $s \in [0, t]$.

Recall that $\Omega^0 = \{M \in \Omega: |[M]| = p^{t-s}n\}$, \cf \autoref{definition-omega}.(1).

\begin{remark without proof}\label{remark-a_0-direkt}
	Recall that $a_0 = |\overline\Omega^0|$, \cf \autoref{definition-omega}.(3).
	
	We have $$a_0 = \Num(s),$$ either by \autoref{lemma-wielandt} or by \autoref{proposition-al}. The congruence $$\Num(s) \equiv_p 1$$ follows from \autoref{thm-sylow-wielandt} or from \autoref{thm-sub-formel-pl}.
	
	So the number $\Num(s)$ of subgroups $U$ of $G$ with $|U| = p^s$ is congruent to $1$ modulo $p$.
	
	This congruence is the statement of the Theorem of Sylow-Frobenius \bfcite{Syl}{Th. II}.
\end{remark without proof}

%% file: text/BA-a_1.tex
Let $s \in [1, t]$.

Recall that $\Omega_1^k = \{M \in \Omega: 1 \in M,\ |[M]| = p^{t-s+k}n\}$ for $k \in [0, 1]$, \cf Definitions \ref{definition-omega-1-l}.(1) and \ref{definition-omega}.(1).

We have the following situation, \cf \autoref{definition-data}.

\begin{center}
	\begin{tikzcd}[column sep=1ex,row sep=5ex]
		&&\Data^{(1)}\arrow[d, "\DataPhi^{(1)}"]
		&=
		&\Data^{(1),0}\arrow[d, "\DataPhi^{(1),0}"]
		&\sqcup
		&\Data^{(1),1}\arrow[d, "\DataPhi^{(1),1}"] \\
		
		\Omega\arrow[d, "\DataRho"]
		&\supseteq
		&\Omega_1^{[0,1]}\arrow[d, "\DataRho_1^{[0,1]}"]
		&=
		&\Omega_1^0\arrow[d, "\DataRho_1^{0}"]
		&\sqcup
		&\Omega_1^1\arrow[d, "\DataRho_1^{1}"] \\
		
		\overline\Omega
		&\supseteq
		&\overline\Omega^{[0,1]}
		&=
		&\overline\Omega^0
		&\sqcup
		&\overline\Omega^1
	\end{tikzcd}
\end{center}

\begin{remark without proof}\label{remark-a_1-direkt}
	Recall that $a_1 = |\overline\Omega^1|$, \cf \autoref{definition-omega}.(3). We have $$a_1 = \tfrac 1p (\TvBin(1)\Num(s-1) - \Num(s-1, s)),$$ \cf \autoref{proposition-al}.
	
	We can also see this by a direct calculation, obtaining
	\begin{align*}
		a_1 &\mystackrel{D. \ref{definition-omega}.(3)}{=}
		|\overline\Omega^1| \\
		&\mystackrel{L. \ref{lemma-fibres-phi-bar}.(2)}{=}
		\tfrac 1p |\Data^{(1),1}| \\
		&\mystackrel{D. \ref{definition-data}.(3)}{=}
		\tfrac 1p (|\Data^{(1)}| - |\Data^{(1),0}|) \\
		&\mystackrel{D. \ref{definition-data}.(1)}{=}
		\tfrac 1p (\TvBin(1)\Num(s-1) - |\Data^{(1),0}|) \\
		&\mystackrel{L. \ref{lemma-bijection-sub-sml}}{=}
		\tfrac 1p (\TvBin(1)\Num(s-1) - \Num(s-1, s)).\qedhere
	\end{align*}
\end{remark without proof}

%% file: text/BA-a_2.tex
Let $s \in [2, t]$.

Recall that $\Omega_1^k = \{M \in \Omega: 1 \in M,\ |[M]| = p^{t-s+k}n\}$ for $k \in [0, 2]$, \cf Definitions \ref{definition-omega-1-l}.(1) and \ref{definition-omega}.(1).

We have the following situation, \cf \autoref{definition-data}.

\begin{center}
	\begin{tikzcd}[column sep=1ex,row sep=5ex]
		&&\Data^{(2)}\arrow[d, "\DataPhi^{(2)}"]
		&=
		&\Data^{(2),0}\arrow[d, "\DataPhi^{(2),0}"]
		&\sqcup
		&\Data^{(2),1}\arrow[d, "\DataPhi^{(2),1}"]
		&\sqcup
		&\Data^{(2),2}\arrow[d, "\DataPhi^{(2),2}"] \\
		
		\Omega\arrow[d, "\DataRho"]
		&\supseteq
		&\Omega_1^{[0,2]}\arrow[d, "\DataRho_1^{[0,2]}"]
		&=
		&\Omega_1^0\arrow[d, "\DataRho_1^{0}"]
		&\sqcup
		&\Omega_1^1\arrow[d, "\DataRho_1^{1}"]
		&\sqcup
		&\Omega_1^2\arrow[d, "\DataRho_1^{2}"] \\
		
		\overline\Omega
		&\supseteq
		&\overline\Omega^{[0,2]}
		&=
		&\overline\Omega^0
		&\sqcup
		&\overline\Omega^1
		&\sqcup
		&\overline\Omega^2
	\end{tikzcd}
\end{center}

\begin{remark without proof}\label{remark-a_2-direkt}
	Recall that $a_2 = |\overline\Omega^2|$, \cf \autoref{definition-omega}.(3).
	We have $$\textstyle a_2 = \tfrac 1{p^2} (\TvBin(2)\Num(s-2) - \Num(s-2, s) - \TvBin(1)\Num(s-2, s-1) + \Num(s-2, s-1, s)),$$ \cf \autoref{proposition-al}.
	
	We can also see this by a direct calculation.
	
	First, note that
	\allowdisplaybreaks
	\[\def\arraystretch{1.5}
		\begin{array}{lcccl}
			|\Data^{(2),1}|
			&\mystackrel{D. \ref{definition-data}.(3)}{=}
			&&\displaystyle\sum_{M \in \Omega_1^1} &|(\DataPhi^{(2),1})^{-1}(\{M\})| \\
			&\mystackrel{L. \ref{lemma-bijection-sub-sml-stab}}{=}
			&&\displaystyle\sum_{M \in \Omega_1^1} &|\Sub(s-2, \Stab(M))| \\
			&\mystackrel{L. \ref{lemma-bijection-omega-1-l}}{=}
			&&\displaystyle\sum_{(U, \{g_2, \dots, g_p\}) \in \Data^{(1),1}} &|\Sub(s-2, U)| \\
			&\mystackrel{D. \ref{definition-data}.(3)}{=}
			&&\displaystyle\sum_{(U, \{g_2, \dots, g_p\}) \in \Data^{(1)}} &|\Sub(s-2, U)| \\
			&\mystackrel{}{\phantom{=}}
			&-&\displaystyle\sum_{(U, \{g_2, \dots, g_p\}) \in \Data^{(1),0}} &|\Sub(s-2, U)| \\
			&\mystackrel{D. \ref{definition-data}.(1), R. \ref{remark-binomial}}{=}
			&&\displaystyle\sum_{U \in \Sub(s-1)} &\TvBin(1)|\Sub(s-2, U)| \\
			&\mystackrel{}{\phantom{=}}
			&-&\displaystyle\sum_{(U, \{g_2, \dots, g_p\}) \in \Data^{(1),0}} &|\Sub(s-2, U)| \\
			&\mystackrel{L. \ref{lemma-bijection-sub-sml}}{=}
			&&\displaystyle\sum_{U \in \Sub(s-1)} &\TvBin(1)|\Sub(s-2, U)| \\
			&\mystackrel{}{\phantom{=}}
			&-&\displaystyle\sum_{(U, V) \in \Sub(s-1, s)} &|\Sub(s-2, U)| \\
%			&\mystackrel{}{=}
%			&&\displaystyle\sum_{U \in \Sub(s-1)} &\TvBin(1)|\Sub(s-2, U)| \\
%			&\mystackrel{}{\phantom{=}}
%			&-&\displaystyle\sum_{U \in \Sub(s-1)} &|\Sub(s-2, U)||\Sub(U, s)| \\
			&\mystackrel{D. \ref{definition-sub-num}.(2)}{=}
			&\rlap{$\TvBin(1)\Num(s-2, s-1) - \Num(s-2, s-1, s),$}
		\end{array}
	\]
	\cf also \autoref{cor-cardinality-d-k-l}. Then
	\begin{align*}
		a_2 &\mystackrel{D. \ref{definition-omega}.(3)}{=}
		|\overline\Omega^2| \\
		&\mystackrel{L. \ref{lemma-fibres-phi-bar}.(2)}{=}
		\tfrac 1{p^2} |\Data^{(2),2}| \\
		&\mystackrel{D. \ref{definition-data}.(3)}{=}
		\tfrac 1{p^2} (|\Data^{(2)}| - |\Data^{(2),0}| - |\Data^{(2),1}|) \\
		&\mystackrel{D. \ref{definition-data}.(1)}{=}
		\tfrac 1{p^2} (\TvBin(2)\Num(s-2) - |\Data^{(2),0}| - |\Data^{(2),1}|) \\
		&\mystackrel{L. \ref{lemma-bijection-sub-sml}}{=}
		\tfrac 1{p^2} (\TvBin(2)\Num(s-2) - \Num(s-2, s) - |\Data^{(2),1}|) \\
		&\mystackrel{}{=}
		\tfrac 1{p^2} (\TvBin(2)\Num(s-2) - \Num(s-2, s) \\
		&\mystackrel{}{\phantom{=}}
		\phantom{\tfrac 1{p^2} (}- (\TvBin(1)\Num(s-2, s-1) - \Num(s-2, s-1, s))). \qedhere
	\end{align*}
\end{remark without proof}